\documentclass[11pt]{article}

\usepackage{amssymb}
\usepackage{latexsym}
\usepackage[tbtags]{amsmath}
\usepackage{epsfig}
\usepackage{amscd}
\usepackage{amsthm}
\usepackage{multirow}
\usepackage{float}
\usepackage{footmisc}
\usepackage{color}
\usepackage{enumerate}
\usepackage{bm}
\usepackage{wrapfig}
\usepackage{setspace}
\usepackage{adjustbox} 
\usepackage{xparse}
\usepackage{tcolorbox}
\usepackage{lipsum}
\usepackage{blindtext}
\usepackage{mathtools}
\usepackage[titletoc,title]{appendix}
\usepackage[T1]{fontenc}
\usepackage{dsfont}
\usepackage{booktabs}
\usepackage{natbib}
\usepackage{algorithm}
\usepackage{enumitem}
\usepackage{algpseudocode}
\usepackage{url}
\usepackage[font=small,labelfont=bf]{caption}
\usepackage{authblk} 
\usepackage{hyperref}[] 
\hypersetup{
colorlinks = true,
linkcolor = blue, 
filecolor = magenta,
urlcolor = blue, 
citecolor = blue 
}

\usepackage{sectsty}

\newcommand{\mP}{\mathbb{P}}
\newcommand{\mE}{\mathbb{E}}

\newcommand{\mV}{\mathrm{Var}}

\newcommand{\tr}{\text{tr}}

\newcommand{\bK}{\boldsymbol{\mathrm{K}}}
\newcommand{\bL}{\boldsymbol{\mathrm{L}}}
\newcommand{\bOne}{\boldsymbol{\mathrm{1}}}

\newcommand{\indep}{\rotatebox[origin=c]{90}{$\models$}} 
\DeclarePairedDelimiter\ceil{\lceil}{\rceil}
\DeclarePairedDelimiter\floor{\lfloor}{\rfloor}

\newcommand{\convP}{\overset{p}{\longrightarrow}}

\newtheorem{theorem}{Theorem}
\newtheorem{lemma}{Lemma}

\newtheorem{proposition}{Proposition}

\newtheorem{remark}{Remark}

\allowdisplaybreaks

\begin{document}

 \begin{center}
 	\textbf{\LARGE Conditional Independence Testing for Discrete Distributions: Beyond $\chi^2$- and $G$-tests}	
 
 	\vspace*{.2in}
 	
 	\begin{author}
 		A
 		Ilmun Kim$^{\dagger}$  \quad Matey Neykov$^{\ddagger}$ \quad Sivaraman Balakrishnan$^{\S}$ \quad Larry Wasserman$^{\S}$
 	\end{author}
 
 	\let\thefootnote\relax\footnotetext{The last three authors are listed \href{https://tinyurl.com/random-order}{randomly}.}
 
 	\vspace*{.1in}
 	
 	\begin{tabular}{c}
 		$^{\dagger}$Department of Statistics and Data Science, Yonsei University\\
 		$^{\ddagger}$Department of Statistics and Data Science, Northwestern University\\
 		$^{\S}$Department of Statistics and Data Science, Carnegie Mellon University
 	\end{tabular}

 	\vspace*{.1in}
 	
 	\today
 	
 	\vspace*{.1in}
 \end{center}

\begin{abstract} 
	This paper is concerned with the problem of conditional independence testing for discrete data. In recent years, researchers have shed new light on this fundamental problem, emphasizing finite-sample optimality. The non-asymptotic viewpoint adapted in these works has led to novel conditional independence tests that enjoy certain optimality under various regimes. Despite their attractive theoretical properties, the considered tests are not necessarily practical, relying on a Poissonization trick and unspecified constants in their critical values. In this work, we attempt to bridge the gap between theory and practice by reproving optimality without Poissonization and calibrating tests using Monte Carlo permutations. Along the way, we also prove that classical asymptotic $\chi^2$- and $G$-tests are notably sub-optimal in a high-dimensional regime, which justifies the demand for new tools. Our theoretical results are complemented by experiments on both simulated and real-world datasets. Accompanying this paper is an \texttt{R} package \texttt{UCI} that implements the proposed tests.
\end{abstract}

\section{Introduction}

Conditional independence (CI) is the backbone of diverse fields in statistics, including graphical models~\citep{de2000new,koller2009probabilistic} and causal inference~\citep{spohn1994properties,pearl2014probabilistic,imbens2015causal}. Among several benefits, this fundamental assumption allows us to simplify the structure of a model, thereby increasing interpretability and reducing computational costs. To justify the use of CI assumption, it is of considerable interest to test whether two random variables $X$ and $Y$ are independent after accounting for the effect of another random variable $Z$. Due to its important role, the problem of CI testing has received much attention in the past decade, resulting in numerous exciting new developments~\citep[e.g.,][]{candes2016panning,shah2020hardness,berrett2020conditional,neykov2021minimax,petersen2021testing, liu2022fast,lundborg2022conditional}. See \cite{li2020nonparametric} and \cite{chatterjee2022survey} for recent reviews. However, most of the recent work is dedicated to continuous data and the importance of discrete CI testing is relatively overlooked.

In discrete settings, two commonly used methods are the $\chi^2$-test~\citep{pearson1900x} and the $G$-test~\citep{mcdonald2014g}, and their asymptotic equivalence is well-known under regularity conditions~\citep[e.g.,~Chapter 14 of][]{bishop2007discrete}. When $X$ and $Y$ are binary, the Cochran--Mantel--Haenszel test~\citep{agresti2003categorical} is another popular method for CI testing. Despite their popularity, these methods are
asymptotic in nature, frequently calibrated by their limiting null distributions. Therefore their finite-sample validity remains questionable. This miscalibration issue becomes more serious in high-dimensional regimes where the number of categories can be significantly larger than the sample size. Besides, the power of these methods is not well-understood except in classical fixed-dimensional settings.

For discrete CI testing, \cite{canonne2018testing} put forward two testing algorithms and analyze their sample complexity from a non-asymptotic perspective. Their sample complexity results are further complemented by matching lower bounds, demonstrating optimality of their procedures in some regimes. In spite of these technical advances, their approach poses several practical challenges. First, their results rely on a Poissonization trick where the sample size is treated as a Poisson random variable. This assumption greatly simplifies the theoretical analysis, but is untenable in practice. Another issue worth highlighting is the dependence of the test on unspecified constants in their critical values. In many statistical applications, the type I error is a greater concern than the type II error. It is therefore desirable to set a critical value in such a way as to maximize the power, while \emph{tightly} controlling the type I error. However, it is unclear from \cite{canonne2018testing} how to modify their tests to meet this criteria, thereby leaving room for improvement from a practical perspective. 
Indeed, this issue was the main motivation of recent work of \cite{kim2022minimax,kim2021local} that advocates the use of permutation methods in two-sample and (both unconditional and conditional) independence testing problems.

With these issues in mind, our work makes the following contributions: (i)~In Theorem~\ref{Theorem: Multinomial Sampling}, we depoissonize the sample complexity results of \cite{canonne2018testing} and establish the same theoretical guarantees under the standard sampling setting. On a technical level, the challenge lies in dealing with the complicated dependence structure of multinomial samples. We overcome this difficulty using the negative association property of multinomial distributions~\citep[][]{joag1983negative}. (ii)~We further make the algorithms of \cite{canonne2018testing} practical by leveraging the permutation method to calibrate test statistics. This resampling approach completely removes the issue arising from unspecified constants, and provably controls the finite-sample type I error. In Theorem~\ref{Theorem: Multinomial Sampling using Permutation}, we prove that Monte Carlo permutation tests achieve the same sample complexity as the theoretical tests of \cite{canonne2018testing}. (iii)~The considered test statistics are linear combinations of fourth order U-statistics, which can be daunting computationally. We address this computational concern by presenting alternative linear time expressions in Proposition~\ref{Proposition: Computationally efficient formula}.~(iv)~We also prove an independent result that demonstrates sub-optimality of asymptotic $\chi^2$- and $G$-tests in their power performance. This negative result naturally inspires efforts to develop new CI tests that perform better than the classical ones. (v) Finally, we provide extensive simulation results that demonstrate the practical value of the proposed methods in Section~\ref{Section: Numerical analysis}, and the algorithms are available in the \texttt{R} package \texttt{UCI}.\footnote{publicly available at github repository:  \href{https://github.com/ilmunk/UCI}{\texttt{https://github.com/ilmunk/UCI}}}


Our work is related to \cite{tsamardinos2010permutation} who warn about the risk of asymptotic calibration for $\chi^2$- and $G$-tests, and further highlight benefits of the permutation procedure in type I error control. The risk of asymptotic calibration has also been discussed in other testing problems, such as those studied in \cite{balakrishnan2018hypothesis,balakrishnan2019hypothesis,kim2022minimax}. In line with this research, we prove the negative result of asymptotic $\chi^2$- and $G$-tests, and demonstrate attractive properties of the permutation method both in type I and II error control. Another related work is \cite{berrett2021usp} where the authors propose a permutation test based on a U-statistic for unconditional independence testing. Concurring with our view, \cite{berrett2021usp} put an emphasis on the permutation approach for practical calibration and demonstrate the competitive performance of their proposal, coined \texttt{USP} test, over $\chi^2$- and $G$-tests. In fact, when the conditional variable is degenerate~(i.e.,~$Z$ takes a single value), one of our practical proposals becomes exactly the same as that of \cite{berrett2021usp}. In this sense, our work can be considered as an extension of \cite{berrett2021usp} to CI testing. We also refer to \cite{agresti1992survey,yao1993exact} that discuss exact inference methods for contingency tables. It is worth pointing out that the current paper builds on our prior work~\citep{kim2021local}, which proves that the sample complexity results of \cite{canonne2018testing} continue to hold using permutation tests. However, the analysis of \cite{kim2021local} relies on Poissonization and also makes use of (computationally expensive) full permutation tests. The current work deviates from \cite{kim2021local} by removing Poissonization and employing a more computationally efficient permutation test via Monte Carlo sampling. We also propose a new permutation test, called \texttt{wUCI}-test, that avoids sample splitting and illustrate its competitive finite sample performance under a variety of settings.

\paragraph{Organization.} The rest of this paper is organized as follows. In Section~\ref{Section: Background}, we set the stage by presenting some background information on sample complexity and Poissonization. Section~\ref{Section: Test statistics} describes the test statistics that we study, and verifies that they can be computed in linear time. Section~\ref{Section: Main theoretical results} contains our main theoretical results including depoissonization and sub-optimality of $\chi^2$- and $G$-tests. In Section~\ref{Section: Numerical analysis}, we demonstrate the empirical performance of the proposed methods based on simulated and real-world datasets, before concluding in Section~\ref{Section: Discussion}. All the proofs of our results are relegated to the Appendix.

\paragraph{Notation.} For a positive integer $a$, we use the shorthand $[a] = \{1,\ldots,a\}$. The conditional independence of $X$ and $Y$ given $Z$ is denoted as $X \, \indep \,Y \,|\, Z$. Given two discrete distributions $p$ and $q$, we write the $L_1$ distance between $p$ and $q$ as $\|p-q\|_1$. We say that random variables $X_1,\ldots,X_n$ are i.i.d. when they are independent and identically distributed. For two positive sequences $a_n$ and $b_n$, we write $a_n \asymp b_n$ if it holds that $C_1 \leq a_n/b_n \leq C_2$ for some positive constants $C_1$ and $C_2$, and for all $n$. We also write $a_n = O(b_n)$ or $a_n \lesssim b_n$ to indicate that $a_n \leq C b_n$ for some positive constant $C$ independent of $n$. 

\section{Background} \label{Section: Background}
Before presenting our main results, we start by building some background knowledge on sample complexity and Poissonization.
\subsection{Setting the stage}
\begin{figure}[t!]
	\begin{center}		
		\includegraphics[width=36em]{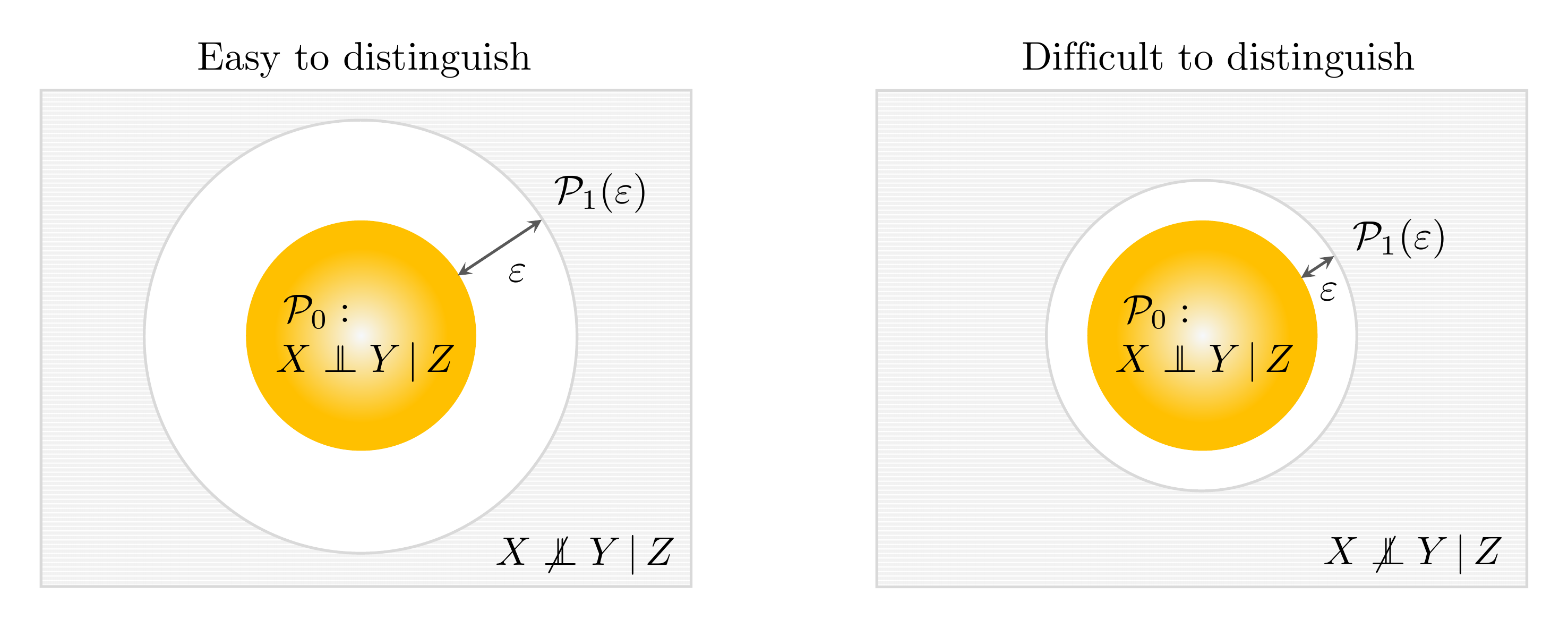} 
		\caption{\small Schematic of the hypotheses of CI testing. The space of null distributions is $\varepsilon$ far from the space of alternative distributions where the distance parameter $\varepsilon$ controls the difficulty of the problem.} \label{Figure: hypothesis}
	\end{center}
\end{figure}

Consider the set of discrete distributions of $(X,Y,Z)$ on a domain $[\ell_1] \times [\ell_2] \times [d]$, denoted by $\mathcal{P}$. Let $\mathcal{P}_0$ be the subset of $\mathcal{P}$ such that $X \, \indep \, Y \,|\, Z$. Given $n$ i.i.d.~random vectors $\{(X_i,Y_i,Z_i)\}_{i=1}^n$ drawn from $p_{X,Y,Z} \in \mathcal{P}$, our goal is to distinguish
\begin{align} \label{Eq: hypothesis}
	H_0: p_{X,Y,Z} \in \mathcal{P}_0 \quad \text{versus} \quad H_1: p_{X,Y,Z} \in \mathcal{P}_1(\varepsilon) = \Big\{p \in \mathcal{P} : \inf_{q \in \mathcal{P}_0} \|p - q\|_1 \geq \varepsilon \Big\},
\end{align}
where $\varepsilon >0$ is a distance parameter (Figure~\ref{Figure: hypothesis} for a pictorial description). When the sample size $n$ is too small, no valid test can reliably differentiate the null from the alternative. On the other hand, when the sample size $n$ is too large, the problem becomes trivial, resulting in many successful tests. A natural question is then to determine the smallest $n$ required to achieve the desired level of testing accuracy for an optimal test. This concept is known as the optimal sample complexity formally defined below.

\paragraph{Optimal sample complexity.}  We define a test $\phi$ which maps from the samples to a binary outcome as $\phi: \{(X_i,Y_i,Z_i)\}_{i=1}^n \mapsto \{0,1\}$. For some fixed $\alpha \in (0,1)$, let $\Phi_{\alpha}$ denote the set of level $\alpha$ tests such that for each $\phi \in \Phi_{\alpha}$,
\begin{align*}
	\sup_{p \in \mathcal{P}_0} \mP_{p}(\phi = 1) \leq \alpha.
\end{align*}
The minimax risk is the worst-case type II error of an optimal level $\alpha$ test defined as 
\begin{align*}
	R_n(\varepsilon,\alpha) = \inf_{\phi \in \Phi_{\alpha}} \sup_{p \in \mathcal{P}_1(\varepsilon)} \mP_{p}(\phi = 0).
\end{align*}
The difficulty of the problem can be characterized as the minimum number of samples that makes the minimax risk bounded by some fixed constant $\beta \in (0,1-\alpha)$\footnote{Throughout this paper, we treat the target type I and II errors $\alpha$ and $\beta$ as universal constants, e.g.,~$\alpha=\beta=0.05$.}. This minimum number of samples is called the optimal sample complexity given as
\begin{align*}
	n^\star = \inf \Big\{ n : R_n(\varepsilon,\alpha) \leq \beta \Big\}.
\end{align*}
We say that a level $\alpha$ test $\phi$ is rate optimal in sample complexity if 
\begin{align*}
	n^\star \asymp \inf\bigg\{ n : \sup_{p \in \mathcal{P}_1(\varepsilon)} \mP_{p}(\phi = 0) \leq \beta \bigg\}.
\end{align*}
In practice, an optimal test whose risk is exactly equal to $R_n(\varepsilon,\alpha)$ is mostly inaccessible. For this reason, we instead aim to find a rate optimal test.

\begin{remark} \leavevmode \normalfont
	\begin{itemize}
		\item (\emph{Finite-sample optimality}) We highlight that the finite-sample minimax optimality considered in this paper is distinct from traditional asymptotic optimality. In classical asymptotic theory, the performance of a test is typically measured against contiguous alternatives, which converge to the null hypothesis at a $\sqrt{n}$-rate. While this type of local asymptotic analysis enables precise power comparisons, it does not provide useful insight when the underlying distribution is outside of regular, fixed-dimensional models. This limitation has been emphasized in \citet{balakrishnan2018hypothesis} and \cite{arias2018remember} that advocate for the finite-sample minimax framework as in our work. This framework provides quantitative information in both low- and high-dimensional cases, offering a valuable alternative to traditional asymptotic approaches.
		\item (\emph{Choice of metrics}) As is well-known, one half of the $L_1$ distance between two distributions $p$ and $q$ is the same as their total variation (TV) distance:
		\begin{align*}
			\frac{1}{2}\|p - q\|_1 = \mathrm{TV}(p,q) = \sup_{A} |p(A) - q(A)|, 
		\end{align*} 
		where the supremum is taken over all possible measurable sets. The optimal sample complexity depend crucially on the choice of metrics in the definition of the alternative hypothesis~(\ref{Eq: hypothesis}). Since our aim is simply to depoissonize the previous results, we focus on the $L_1$ distance or equivalently the TV distance throughout this paper. We refer the reader to the recent work of \cite{neykov2023nearly} that investigates the optimal sample complexity for CI testing in terms of the Wasserstein distance. 
	\end{itemize}
\end{remark}

\subsection{Poissonization}

Poissonization is now a standard technique to study the sample complexity in hypothesis testing~\citep[e.g.,][]{valiant2011testing,chan2014optimal,valiant2017automatic,balakrishnan2019hypothesis}. The idea itself is old in statistics and probability theory, dating back at least to \cite{kac1949deviations}. See also Chapter 3.5 of \cite{vaart1996weak}. It is used in theoretical works on hypothesis testing as a trick to simplify several calculations in dealing with categorical data. In particular, it is well-known that when the data are generated by Poisson sampling (Algorithm~\ref{Algorithm: Poisson sampling}), the number of samples falling into disjoint sets are mutually independent. This independence property lies at the heart of deriving various existing results of the sample complexity in distribution testing. See \cite{canonne2022topics} for a recent review.

\begin{algorithm} \caption{Poisson Sampling} \label{Algorithm: Poisson sampling}
	\textbf{Input:} For a fixed $n \in \mathbb{N}$ and a distribution $p_{X,Y,Z}$
	\begin{enumerate}
		\item Draw $N \sim \mathrm{Poisson}(n)$.
		\item Generate i.i.d.~$(X_1,Y_1,Z_1),\ldots,(X_N,Y_N,Z_N)$ random variables from $p_{X,Y,Z}$.
	\end{enumerate}
	\textbf{Return:} $\{(X_i,Y_i,Z_i)\}_{i=1}^N$
\end{algorithm}

\subsection{General recipe and related issue} \label{Section: A general recipe} 
If constant factors are not the main concern, there are straightforward ways of transferring the sample complexity obtained from Poisson sampling to the usual sampling scenario with a fixed sample size. One concrete procedure, described in \cite{neykov2021minimax}, is as follows. Given $n$ i.i.d.~copies of $(X,Y,Z)$, 
\begin{enumerate}
	\item Draw $N \sim \mathrm{Poisson}\bigl(\frac{n}{2}\bigr)$.
	\item If $N > n$, then accept the null hypothesis. 
	\item If $N \leq n$, perform a test based on $(X_1,Y_1,Z_1),\ldots,(X_N,Y_N,Z_N)$ and return a result. 
\end{enumerate}
Let $\phi$ be a generic test function using $(X_1,Y_1,Z_1),\ldots,(X_N,Y_N,Z_N)$ and then the resulting test from the above procedure can be concretely written as $\phi^\ast = \mathds{1}(N \leq n) \phi$. Suppose that the test $\phi$ has the type I error as well as the type II error bounded by $\alpha$ and $\beta$, respectively. Then the corresponding test $\phi^\ast$ has the type I and II error bounds as
\begin{align*}
	\sup_{p \in \mathcal{P}_0} \mE_{p}[\phi^\ast] \leq \alpha \quad \text{and} \quad \sup_{p \in \mathcal{P}_1(\varepsilon)} \mE_{p}[1 - \phi^\ast] \leq  \beta + \mP(N > n).
\end{align*}
Since the Poisson distribution is tightly concentrated around its mean, the additional term in the type II error of $\phi^\ast$ can be made small when $n$ is relatively large compared to $\beta$.\footnote{More precisely, an exponential tail bound for a Poisson random variable ensures that $\mP(N > n) \leq e^{-n/8}$~\citep[e.g.,~Theorem A.8 of][]{canonne2022topics} so that if $n \geq 8 \log(1/\beta)$, the type II error is bounded above by $2\beta$.} Therefore one can transfer the sample complexity obtained under Poissonization to the usual sampling scheme up to a constant factor. 

Nevertheless, due to an inefficient use of the data as well as a non-trivial chance of accepting the null irrelevant to the data, practitioners may not necessarily follow this general recipe. To address this concern, there has been an effort to depoissonize the results of distribution testing under Poisson sampling. For instance, \cite{kim2020multinomial} revisits the truncated $\chi^2$ test for goodness-of-fit testing proposed by \cite{balakrishnan2019hypothesis} and establishes the same optimality without Poissonization. We also refer to \cite{jacquet1998analytical} and Chapter 2.5 of \cite{penrose2003random} that present useful depoissonization tools, but mostly for asymptotic results. In this work, we analyze the tests proposed by \cite{canonne2018testing} under the usual i.i.d.~sampling model with a fixed sample size and show that the same sample complexity can be achieved without any further assumption. At the heart of our technique is the negative association of multinomial distributions~\citep{joag1983negative}, which would be potentially useful for depoissonizing other sample complexity results.

\vskip 1em

\noindent \textbf{Multinomial sampling.} To fix terminology, we refer to the usual sampling with fixed sample size $n$ as \emph{multinomial sampling}  in what follows.

\section{Test statistics} \label{Section: Test statistics}
To describe our main results, we first need to recall the test statistics introduced by \cite{canonne2018testing}. As noted by \cite{neykov2021minimax}, their test statistics can be viewed as linear combinations of U-statistics constructed using the observations of $(X,Y)$ in the same category of $Z$. We describe two kinds of U-statistics considered in \cite{canonne2018testing}. 

\subsection{Unweighted U-statistic}
Suppose that we observe $\sigma \geq 4$ i.i.d.~observations of $(X,Y)$ supported on $[\ell_1] \times [\ell_2]$. We let $p_{X,Y}$ denote the joint discrete distribution of $(X,Y)$ and let $p_X$ and $p_Y$ denote the marginal distribution of $X$ and $Y$, respectively. The first U-statistic is an unbiased estimator of the squared $L_2$ distance between $p_{XY}$ and $p_Xp_Y$. In more detail, borrowing the notation from \cite{neykov2021minimax}, let 
\begin{align*}
	\phi_{ij}(qr) = \mathds{1}(X_i = q)\mathds{1}(Y_i = r) - \mathds{1}(X_i= q) \mathds{1}(Y_j=r). 
\end{align*}
For four distinct observations indexed by $i,j,k,l \in [\sigma]$, the kernel of the unweighted U-statistic is defined as
\begin{align*}
	h_{ijkl} = \frac{1}{4!} \sum_{(\pi_1,\pi_2,\pi_3,\pi_4) \in \Pi} \sum_{q \in [\ell_1],r \in [\ell_2]} \phi_{\pi_1 \pi_2}(qr) \phi_{\pi_3 \pi_4} (qr),
\end{align*}
where $\Pi$ is the set of all possible permutations of $(i,j,k,l)$. By the linearity of expectations, we see that $h_{ijkl}$ is an unbiased estimator of the squared $L_2$ distance between $p_{XY}$ and $p_Xp_Y$. Given this kernel and the dataset $D = \{(X_1,Y_1),\ldots,(X_\sigma,Y_\sigma)\}$, the unweighted U-statistic is computed as
\begin{align} \label{Eq: unweighted U}
	U(D) = \binom{\sigma}{4}^{-1} \sum_{i <j <k < l: (i,j,k,l) \in [\sigma]} h_{ijkl}.
\end{align}
It is worth pointing out that $U(D)$ is equivalent to the U-statistic for unconditional independence testing considered in \cite{berrett2021optimal,berrett2021usp,kim2022minimax}. Denote the datasets $D_m =\{(X_i,Y_i) : Z_i = m, i \in [n]\}$ and write the sample size within $D_m$ by $\sigma_m$. The final test statistic for CI testing based on $U(D)$ is then constructed as 
\begin{align*}
	T = \sum_{m \in [d]} \mathds{1}(\sigma_m \geq 4) \sigma_m U(D_m).
\end{align*}	 

\subsection{Weighted U-statistic} \label{Seciton: Weighted U-statistic}
The previous unweighted U-statistic may suffer from high variance especially when the $L_2$ norms of $p_{XY}$ and $p_Xp_Y$ are large, resulting in sub-optimal performance (see the simulation result in Figure~\ref{Figure: Power I} under Scenario 1). To address this issue, \cite{canonne2018testing} employ the flattening idea proposed by \cite{diakonikolas2016new}. This approach involves partitioning heavier bins into multiple smaller pieces in a way to reduce the $L_2$ norms of the transformed distributions, leading to a smaller variance of the test statistic. As noted in \cite{balakrishnan2018hypothesis} and further elaborated by \cite{neykov2021minimax}, the flattening procedure is equivalent to using a carefully designed weighted kernel for a U-statistic.

To describe the weighted U-statistic, assume that the sample size $\sigma = 4 + 4t$ for some $t \in \mathbb{N}$, and let $t_1 = \min(t,\ell_1)$ and $t_2 = \min(t,\ell_2)$. The construction of the weighted U-statistic involves sample splitting where the dataset $D = \{(X_1,Y_1),\ldots,(X_\sigma,Y_\sigma)\}$ is split into $D_X = \{X_i: i \in [t_1]\}$, $D_Y = \{Y_i: t_1 + 1 \leq i \leq t_1 + t_2\}$ and $D_{X,Y} = \{(X_i,Y_i): 2t + 1 \leq i \leq \sigma\}$ of sizes $t_1$, $t_2$ and $4 + 2t$, respectively. Define $a_q$ (and $a_r'$) as the number of observations equal to $q$ (and $r$) in $D_X$ (and $D_Y$). For four distinct observations indexed by $i,j,k,l$ in $D_{X,Y}$, consider a weighted kernel given as
\begin{align*}
	h_{ijkl}^{\boldsymbol{a}} = \frac{1}{4!} \sum_{(\pi_1,\pi_2,\pi_3,\pi_4) \in \Pi} \sum_{q \in [\ell_1],r \in [\ell_2]} \frac{\phi_{\pi_1 \pi_2}(qr) \phi_{\pi_3 \pi_4} (qr)}{(1+a_q)(1+a_r')}.
\end{align*}
Due to the independence among the split datasets, the conditional expectation of $h_{ijkl}^{\boldsymbol{a}}$ taken over $D_{X,Y}$ is the square of the $L_2$ distance between $p_{XY}$ and $p_Xp_Y$ weighted by $(1+a_q)(1+a_r')$. Given the kernel, the weighted U-statistic is computed as 
\begin{align} \label{Eq: weighted U}
	U_W^{\boldsymbol{a}}(D) = \binom{2t+4}{4}^{-1} \sum_{i <j <k < l: (i,j,k,l) \in D_{X,Y}} h_{ijkl}^{\boldsymbol{a}}, 
\end{align}
where $(i,j,k,l) \in D_{X,Y}$ indicates taking four observations indexed by $(i,j,k,l)$ from the dataset $D_{X,Y}$. As before, denote the datasets $D_m =\{(X_i,Y_i) : Z_i = m, i \in [n]\}$ and write the sample size within $D_m$ by $\sigma_m$. By further writing $\omega_m = \sqrt{\min(\sigma_m,\ell_1) \min(\sigma_m,\ell_2)}$, the final test statistic for CI testing based on $U_W^{\boldsymbol{a}}(D)$ is then constructed as 
\begin{align*}
	T_W= \sum_{m \in [d]} \mathds{1}(\sigma_m \geq 4) \sigma_m \omega_m U_W^{\boldsymbol{a}}(D_m).
\end{align*}

\paragraph{Practical considerations.} As shown in \cite{canonne2018testing}, the test based on $T_W$ achieves the optimal sample complexity in broader regimes than that based on $T$. Despite its attractive theoretical properties, the resulting test may experience a loss of practical power due to an inefficient use of the data arising from sample splitting. To mitigate this issue, we introduce another weighted U-statistic without sample splitting. First, it is worth trying to understand the weight $(1+a_q)(1+a_r')$ at a population level. In particular, the conditional expectation of the weight given the random sample size $t$ is $\bigl(1 + \min\{t, \ell_1\}p_X(q)\bigr) \bigl(1 + \min\{t, \ell_1\} p_Y(r)\bigr)$. This motivates us to consider a weight kernel 
	\begin{align*}
	h_{ijkl}^{\boldsymbol{b}} = \frac{1}{4!} \sum_{(\pi_1,\pi_2,\pi_3,\pi_4) \in \Pi} \sum_{q \in [\ell_1],r \in [\ell_2]} \frac{\phi_{\pi_1 \pi_2}(qr) \phi_{\pi_3 \pi_4} (qr)}{(1+b_q)(1+b_r')},
\end{align*}
where $b_q = \min\{\sigma, \ell_1\} \sigma^{-1} \sum_{i=1}^\sigma \mathds{1}(X_i = q)$ and $b_r' = \min\{\sigma, \ell_2\} \sigma^{-1} \sum_{i=1}^\sigma \mathds{1}(Y_i = r)$. The resulting weighted U-statistic is then computed as 
\begin{align}   \label{Eq: weighted U without splitting}
	U_W^{\boldsymbol{b}} = \binom{\sigma}{4}^{-1} \sum_{i <j <k < l: (i,j,k,l) \in [\sigma]} h_{ijkl}^{\boldsymbol{b}}.
\end{align}
Similarly as before, the final CI test statistic based on $U_W^{\boldsymbol{b}}(D)$ is given as 
\begin{align*}
	T_W^\dagger = \sum_{m \in [d]} \mathds{1}(\sigma_m \geq 4) \sigma_m \omega_m U_W^{\boldsymbol{b}}(D_m).
\end{align*}
We coin the permutation test based on $T_W^\dagger$ as \texttt{wUCI}-test, and formally introduce the procedure in Algorithm~\ref{Algorithm: Weighted U-statistic without splitting} of Appendix~\ref{Section: Algorithm}. 

\begin{remark}[Connection with the truncated $\chi^2$-test] \normalfont \label{Remark: discussion on weights}
	As pointed out by several authors \citep{haberman1988warning,balakrishnan2019hypothesis,kim2020multinomial}, the classical $\chi^2$-test for goodness-of-fit testing can easily break down for sparse multinomial data. To address this problem, \cite{balakrishnan2019hypothesis} introduce a modification of the $\chi^2$-test by using a truncated weight function and prove its minimax optimality. Interestingly, the weight $(1+b_q)(1+b_r')$ that we consider is closely connected to the truncated weight of \cite{balakrishnan2019hypothesis} and may be regarded as an empirical counterpart for independence testing. To explain, let us simply focus on the first term in the product weight and notice that
	\begin{align*}
		\ell_1^{-1}(1+b_q) = \frac{1}{\ell_1}+ \min\biggl( \frac{\sigma}{\ell_1},1\biggr) \widehat{p}_X(q) \asymp \max \biggl\{ \frac{1}{\ell_1}, \min\biggl( \frac{\sigma}{\ell_1},1\biggr) \widehat{p}_X(q) \biggr\}, 
	\end{align*}
	where $\widehat{p}_X(q) =  \sigma^{-1} \sum_{i=1}^\sigma \mathds{1}(X_i = q)$. When $\sigma \ell_1^{-1}$ is large, the right-hand side of the above equation approximates $\max\{\ell_1^{-1}, p_X(q)\}$, which is exactly the same as the truncated weight in \cite{balakrishnan2019hypothesis} for goodness-of-fit testing.
\end{remark}

\subsection{Linear time expression} \label{Section: Linear time expression}
The original forms of the aforementioned U-statistics take $O(\sigma^4 \ell_1 \ell_2)$ time to compute, which can be prohibitive for large $\sigma, \ell_1, \ell_2$. Luckily, this computational complexity can be reduced to $O(\sigma)$ by exploiting a contingency table representation. A computationally convenient form of the unweighted U-statistic $U(D)$ is already given by \cite{canonne2018testing,berrett2021optimal,berrett2021usp}. We also note that an alternative form of $U_W^{\boldsymbol{a}}(D)$ is provided in \cite{canonne2018testing}, but a naive calculation of their expression takes at least $O(\sigma \ell_1^2\ell_2^2)$ time. Here we present a general expression for the U-statistics and explain that it can be run in linear time \emph{independent of} $\ell_1$ and $\ell_2$. 

To this end, let us set some notation. Let $\boldsymbol{\eta} = \{\eta_{1},\ldots,\eta_{\ell_1}\}$ and $\boldsymbol{\upsilon} = \{\upsilon_{1},\ldots,\upsilon_{\ell_2}\}$ be some weight vectors with non-zero components. Given the dataset $D = \{(X_1,Y_1),\ldots,(X_\sigma,Y_\sigma)\}$, consider four distinct observations indexed by $i,j,k,l \in [\sigma]$ and define the weighted kernel associated with $\boldsymbol{\eta}$ and $\boldsymbol{\upsilon}$ as 
\begin{align*}
	h_{ijkl}^{\boldsymbol{\eta},\boldsymbol{\upsilon}} = \frac{1}{4!} \sum_{(\pi_1,\pi_2,\pi_3,\pi_4) \in \Pi} \sum_{q \in [\ell_1], r \in [\ell_2]} \frac{\phi_{\pi_1 \pi_2}(qr) \phi_{\pi_3 \pi_4} (qr)}{\eta_q \upsilon_r}.
\end{align*}
It is clear that the above kernel is equivalent to $h_{ijkl}$ when $\boldsymbol{\eta}$ and $\boldsymbol{\upsilon}$ are vectors with each entry equal to one. Similarly, when $\boldsymbol{\eta}$ and $\boldsymbol{\upsilon}$ are defined with $1+a_q$ and $1+a_r'$, respectively, then the above kernel corresponds to $h_{ijkl}^{\boldsymbol{a}}$. The U-statistic based on $h_{ijkl}^{\boldsymbol{\eta},\boldsymbol{\upsilon}}$ is denoted by $U^{\boldsymbol{\eta},\boldsymbol{\upsilon}}_W (D)$, which is similarly computed as in (\ref{Eq: unweighted U}). For $q \in [\ell_1]$ and $r \in [\ell_2]$, define $o_{qr} =\sum_{i=1}^n \mathds{1}(X_i=q)\mathds{1}(Y_i=r)$, $o_{q+} = \sum_{r=1}^{\ell_2} o_{qr}$ and $o_{+r} = \sum_{q=1}^{\ell_1} o_{qr}$. With this notation in place, we give an alternative expression of $U^{\boldsymbol{\eta},\boldsymbol{\upsilon}}_W (D)$ as follows.

\begin{proposition}[Alternative expression] \label{Proposition: Computationally efficient formula}
	The U-statistic $U^{\boldsymbol{\eta},\boldsymbol{\upsilon}}_W(D)$ can be written as 
	\begin{align*}
		U^{\boldsymbol{\eta},\boldsymbol{\upsilon}}_W (D) ~=~ \frac{1}{\sigma(\sigma-3)}\Biggl[ A_1 + \frac{1}{(\sigma-1)(\sigma-2)} A_2 - \frac{2}{\sigma-2}A_3 \Biggr]
		\end{align*}
		where
		\begin{align*}
			& A_1 =  \sum_{q=1}^{\ell_1}  \sum_{r=1}^{\ell_2} \Biggl( \frac{o_{qr}^2 - o_{qr}}{\eta_q \upsilon_r} \Biggr), \quad A_2 =  \sum_{q=1}^{\ell_1} \Biggl( \frac{o_{q +}^2 - o_{q +}}{\eta_q} \Biggr) \cdot \sum_{r=1}^{\ell_2} \Biggl( \frac{o_{+ r}^2 - o_{+ r}}{\upsilon_r} \Biggr), \\[.5em]
			& A_3 =  \sum_{q=1}^{\ell_1} \sum_{r=1}^{\ell_2} \frac{o_{qr}(o_{q+}o_{+r} - o_{q+} - o_{+r} + 1)}{\eta_q\upsilon_r}.
		\end{align*}
\end{proposition}

We now discuss the average time complexity of $U^{\boldsymbol{\eta},\boldsymbol{\upsilon}}_W (D)$ by assuming that the weight vectors $\boldsymbol{\eta},\boldsymbol{\upsilon}$ are given in advance. First of all, we note that the $\ell_1 \times \ell_2$ contingency table of $D$ is sparse in a sense that it has at most $\sigma$ non-zero entries. Importantly, the zero entries do not affect the calculation of $A_1,A_2, A_3$. Hence we only focus on the non-zero entries of the contingency table, which can be computed in linear time, for instance, by using hash tables~\citep[e.g.,][]{cormen2022introduction}. Similarly, the non-zero row sums and the non-zero column sums of the contingency table can be computed in linear time independent of $\ell_1$ and $\ell_2$. Given the non-zero entries of $\{o_{qr}, o_{q+}, o_{+r}: q \in [\ell_1], r \in [\ell_2]\}$, the computational complexity of the terms $A_1,A_2,A_3$ is linear as their non-zero summands are at most $O(\sigma)$. 

For the unweighted U-statistic $U(D)$, there is no additional cost for computing $\boldsymbol{\eta},\boldsymbol{\upsilon}$ as they are vectors with each entry equal to one. For the weighted U-statistics $U_W^{\boldsymbol{a}}(D)$ and $U_W^{\boldsymbol{b}}(D)$, the weight vectors $\boldsymbol{\eta},\boldsymbol{\upsilon}$ are functions of $o_{q+}$ and $o_{+r}$ (computed on a separate dataset for $U_W^{\boldsymbol{a}}(D)$) and they only require an additional $O(\sigma)$ time to compute. Thus the overall time complexity is still linear.

\section{Theoretical results} \label{Section: Main theoretical results}
Having introduced the test statistics, we are now ready to provide the main theoretical results of this paper. In Section~\ref{Section: Sample complexity without Poissonization}, we establish the same sample complexity of the tests using $T$ and $T_W$ under multinomial sampling. In Section~\ref{Section: Local permutations}, we provide and analyze more practical tests based on permutation procedures.

\subsection{Sample complexity without Poissonization} \label{Section: Sample complexity without Poissonization}
In this subsection, we revisit two main results of \cite{canonne2018testing}, namely Theorem 1.1 and Theorem 1.3 concerning with the sample complexity of a test using $T$ and $T_W$, respectively. 

\vskip 1em

\paragraph{Sample complexity of a test based on $T$:} Suppose that the test statistic $T$ is constructed using $N$ i.i.d.~samples from $p_{X,Y,Z}$ where $N \sim \mathrm{Poisson}(n)$. We reject the null of CI when $T > \zeta \sqrt{\min(n,d)}$ for a sufficiently large constant $\zeta >0$. Then for $\ell_1 = \ell_2 = 2$, Theorem 1.1 of \cite{canonne2018testing} proves that the resulting test has the sample complexity
\begin{align} \label{Eq: sample complexity}
	O\Biggl( \max \Bigg\{ \frac{d^{1/2}}{\varepsilon^2}, \ \min \bigg\{ \frac{d^{7/8}}{\varepsilon}, \ \frac{d^{6/7}}{\varepsilon^{8/7}} \bigg\} \Bigg\} \Biggr).
\end{align}
They also prove that this sample complexity is rate optimal by presenting a matching lower bound. 

\paragraph{Sample complexity of a test based on $T_W$:} The test based on $T$ is not necessarily optimal in the high-dimensional regime where $\ell_1$ and $\ell_2$ can vary with other parameters. As shown in Theorem 1.3 of \cite{canonne2018testing}, a more general result of the sample complexity can be derived by using $T_W$. In particular, given $N$ i.i.d.~samples from $p_{X,Y,Z}$ where $N \sim \mathrm{Poisson}(n)$, we reject the null of CI when $T_W> \zeta' \sqrt{\min(n,d)}$ for a sufficiently large constant $\zeta' > 0$. The sample complexity of the resulting test satisfies
\begin{align} \label{Eq: general sample complexity}
		O\Biggl( \max \Biggl\{  \min \biggl\{ \frac{d^{7/8}\ell_1^{1/4}\ell_2^{1/4} }{\varepsilon}, \ \frac{d^{6/7} \ell_1^{2/7}\ell_2^{2/7} }{\varepsilon^{8/7}}  \biggr\}, \  \frac{d^{3/4}\ell_1^{1/2}\ell_2^{1/2} }{\varepsilon},  \ \frac{d^{2/3}\ell_1^{2/3}\ell_2^{1/3} }{\varepsilon^{4/3}}, \ \frac{d^{1/2}\ell_1^{1/2}\ell_2^{1/2} }{\varepsilon^{2}}  \Biggr\} \Biggr).
\end{align}
Moreover, this upper bound is shown to be optimal, up to constant factors, in a number of regimes. See \cite{canonne2018testing} for a discussion. We note in passing that there was an error in the original proof of Theorem 1.3 in \cite{canonne2018testing}, which has been recently fixed by \cite{kim2022comments}. 

\vskip 1em 

We now depoissonize the previous results and establish the same sample complexity under multinomial sampling.

\begin{theorem}[Multinomial sampling] \label{Theorem: Multinomial Sampling}
	Suppose that we observe $D_n = \{(X_1,Y_1,Z_1),\ldots,(X_n,Y_n,Z_n)\}$ i.i.d.~samples from $p_{X,Y,Z}$ with nonrandom sample size $n$. Then the following two statements hold:
	\begin{enumerate}
		\item Compute $T$ based on $D_n$ and reject the null if $T > \zeta \sqrt{\min(n,d)}$ for a sufficiently large constant $\zeta >0$. Then the resulting test has the sample complexity as in \eqref{Eq: sample complexity} when $\ell_1, \ell_2 = O(1)$.
		\item Compute $T_W$ based on $D_n$ and reject the null if $T_W> \zeta' \sqrt{\min(n,d)}$ for a sufficiently large constant $\zeta' >0$. Then the resulting test has the sample complexity as in \eqref{Eq: general sample complexity}. 
	\end{enumerate}
\end{theorem}

\noindent A few remarks are in order. 

\begin{remark} \leavevmode \normalfont
	\begin{itemize}
		\item The above theorem essentially says that the tests based on $T$ and $T_W$ have the same performance in sample complexity under Poisson sampling and multinomial sampling. This result may not come as a surprise given that a Poisson random variable is strongly concentrated around its mean. See empirical evidence in \cite{kim2021local}. However, the proof under multinomial sampling turns out to be non-trivial, requiring a careful analysis. 
		\item One of the main technical hurdles in the proof is to overcome a lack of independence between random variables in different bins when the sample size is no longer Poisson. The independence property is useful in analyzing the variance of the sum of U-statistics as it leads to zero covariance terms. We address the lack of independence by employing the negative association (NA) property of multinomial random vectors \citep{joag1983negative}. This NA property ensures that the covariance terms are non-positive, which turns out to be enough to ensure the same theoretical guarantees hold under multinomial sampling.
		\item Even though we remove Poissonization, the resulting tests are not necessarily practical. In particular, their critical values depend on unspecified constants $\zeta$ and $\zeta'$. The choice of these constants, resulting in tight control of the type I error, is challenging in practice. We take a further step to address this issue via the permutation method in Section~\ref{Section: Local permutations}, and demonstrate their empirical performance in Section~\ref{Section: Numerical analysis}.
	\end{itemize}
\end{remark}

\subsection{Calibration via permutations} \label{Section: Local permutations}

As mentioned earlier, the tests used in Theorem~\ref{Theorem: Multinomial Sampling} depend on unspecified constants, which raises the issue of practicality. 
This section attempts to address this problem by presenting more practical tests calibrated by the permutation method, and examine their sample complexity under multinomial sampling. In particular, we prove that their sample complexity remains the same as the corresponding (theoretical) tests in Theorem~\ref{Theorem: Multinomial Sampling}. As briefly mentioned earlier, a similar result was established in Theorem 5 of \cite{kim2021local} but under Poisson sampling. In contrast, we do not assume that the sample size follows a Poisson distribution and therefore reduce the gap between theory and practice. We start by describing the testing procedures that we analyze. 

\paragraph{Permutation test using $T$.} This first test compares the test statistic $T$ with its permutation correspondences, and rejects the null if the resulting permutation $p$-value is less than or equal to significance level $\alpha$. To further explain, let $\Pi_{\sigma_m}$ denote the set of all permutations of $[\sigma_m]$ for each $m \in [d]$. Given $\pi$ drawn from $\Pi_{\sigma_m}$, we define $D_m^\pi$ by rearranging $Y$ values in $D_m$ according to $\pi$. More specifically, suppose that we have $D_m = \{(X_1,Y_1),\ldots,(X_{\sigma_m},Y_{\sigma_m})\}$. Then the corresponding permuted dataset becomes $D_m^\pi = \{(X_1,Y_{\pi_1}),\ldots,(X_{\sigma_m},Y_{\pi_{\sigma_m}})\}$. Equipped with this notation, we implement Algorithm~\ref{Algorithm: Unweighted U-statistic} and make a decision based on the output.

\begin{algorithm}[t] \caption{\texttt{UCI}: U-statistic permutation CI test} \label{Algorithm: Unweighted U-statistic}
	\textbf{Input:} Sample~$\{(X_i,Y_i,Z_i)\}_{i=1}^n$, the number of permutations~$B$, significance level $\alpha$
	\vskip .5em
	\textbf{For} $j \in [B]$ \textbf{do}
	\vskip .5em
	\hskip 2em \textbf{For} $m \in [d]$ \textbf{do}
	\begin{itemize}[leftmargin=1.8cm]
		\item Generate $\pi \sim \mathrm{Uniform}(\Pi_{\sigma_m})$ independent of everything else.
		\item Compute $U(D_m^\pi)$ as in (\ref{Eq: unweighted U}) based on the permuted dataset $D_m^\pi$. 
	\end{itemize}
	\hskip 2em \textbf{End}
	\vskip .5em
	\hskip 2em Set $T_j \leftarrow \sum_{m \in [d]} \mathds{1}(\sigma_m \geq 4) \sigma_m U(D_m^\pi)$.
	\vskip .5em
	\textbf{End}
	\begin{itemize}
		\item Set $T \leftarrow \sum_{m \in [d]} \mathds{1}(\sigma_m \geq 4) \sigma_m U(D_m)$ computed without permutations.
		\item Compute the permutation $p$-value
		\begin{align*}
			p_{\mathrm{perm}} = \frac{1}{B+1} \Biggl[ \sum_{j=1}^B \mathds{1}(T_j \geq T) + 1 \Biggr].
		\end{align*}
	\end{itemize}
	\vskip .5em
	\textbf{Output:} Reject $H_0$ if $p_{\mathrm{perm}} \leq \alpha$; otherwise, accept $H_0$. 
\end{algorithm}

\paragraph{Permutation test using $T_W$.} The second test that we analyze calculates its $p$-value by comparing $T_W$ with its permutation correspondences. The overall procedure is similar to the previous one except that it utilizes the half-permutation procedure for a technical reason~\citep[see Remark~6][]{kim2021local}. To explain the procedure, we decompose $D_m$ into $D_{X,m}$, $D_{Y,m}$ and $D_{X,Y,m}$ of size $t_{1,m}$, $t_{2,m}$ and $2t_m+4$, respectively, as in Section~\ref{Seciton: Weighted U-statistic}. Unlike Algorithm~\ref{Algorithm: Unweighted U-statistic}, we only permute $Y$ values within $D_{X,Y,m}$ for each $m \in [d]$ and then evaluate the significance of $T_W$. A more detailed procedure is described in Algorithm~\ref{Algorithm: Weighted U-statistic}.

\begin{algorithm}[t!] \caption{\texttt{wUCI$\_$split}: weighted U-statistic permutation CI test using sample splitting} \label{Algorithm: Weighted U-statistic}
	\textbf{Input:} Sample~$\{(X_i,Y_i,Z_i)\}_{i=1}^n$, the number of permutations~$B$, significance level $\alpha$
	\vskip .5em
	\textbf{For} $j \in [B]$ \textbf{do}
	\vskip .5em
	\hskip 2em \textbf{For} $m \in [d]$ \textbf{do}
	\begin{itemize}[leftmargin=1.8cm]
		\item Generate $\pi \sim \mathrm{Uniform}(\Pi_{2t_m+4})$ independent of everything else.
		\item Define $D_m^\pi =  D_{X,m} \cup D_{Y,m} \cup D_{X,Y,m}^\pi$.
		\item Compute $U_W^{\boldsymbol{a}}(D_m^\pi)$ as in (\ref{Eq: weighted U}) based on the permuted dataset $D_m^\pi$. 
	\end{itemize}
	\hskip 2em \textbf{End}
	\vskip .5em
	\hskip 2em Set $T_{j,W} \leftarrow \sum_{m \in [d]} \mathds{1}(\sigma_m \geq 4) \sigma_m \omega_m U_W^{\boldsymbol{a}}(D_m^\pi)$.
	\vskip .5em
	\textbf{End}
	\begin{itemize}
		\item Set $T_W \leftarrow \sum_{m \in [d]} \mathds{1}(\sigma_m \geq 4) \sigma_m \omega_m U_W^{\boldsymbol{a}}(D_m)$ computed without permutations.
		\item Compute the permutation $p$-value
		\begin{align*}
			p_{\mathrm{perm}} = \frac{1}{B+1} \Biggl[ \sum_{j=1}^B \mathds{1}(T_{j,W} \geq T_W) + 1 \Biggr].
		\end{align*}
	\end{itemize}
	\vskip .5em
	\textbf{Output:} Reject $H_0$ if $p_{\mathrm{perm}} \leq \alpha$; otherwise, accept $H_0$. 
\end{algorithm}
	
\vskip 1em 

Having outlined the permutation procedures, we now discuss their sample complexity. First of all, it is noteworthy that both permutation tests are exact level $\alpha$ in any finite sample scenarios. This simply follows by the fact that the original test statistic and their permutation correspondences are exchangeable under the null. Using this observation, it can be seen that that the resulting $p$-values in Algorithm~\ref{Algorithm: Unweighted U-statistic} and \ref{Algorithm: Weighted U-statistic} are super-uniform~\citep[e.g.,~Lemma 1 of][]{romano2005exact}. The next theorem turns to the type II error and establishes their sample complexity. 

\begin{theorem}[Permutation tests] \label{Theorem: Multinomial Sampling using Permutation}
	Suppose that we observe $D_n = \{(X_1,Y_1,Z_1),\ldots,(X_n,Y_n,Z_n)\}$ i.i.d.~samples from $p_{X,Y,Z}$ with nonrandom sample size $n$. We also assume that the number of random permutations $B$ satisfies $B \geq \max\{4(1-\alpha)\alpha^{-1}, \, 8\alpha^{-2}\log(4\beta^{-1})\}$ where $\alpha$ and $\beta$ are pre-specified type I and II errors, respectively. Then the following two statements hold:
	\begin{enumerate}
		\item The test from Algorithm~\ref{Algorithm: Unweighted U-statistic} has the sample complexity as in \eqref{Eq: sample complexity} when $\ell_1, \ell_2 = O(1)$.
		\item The test from Algorithm~\ref{Algorithm: Weighted U-statistic} has the sample complexity as in \eqref{Eq: general sample complexity}. 
	\end{enumerate}
\end{theorem}
We highlight that the above theorem studies the random permutation tests where $B$ is not required to increase with the sample size $n$. This is in contrast to full permutation tests considered in \cite{kim2021local} that enumerate all possible permutations. We also note that the constant factors in the condition on $B$ is not tight and can be improved with more effort. 

We next turn our attention to $\chi^2$- and $G$-tests and discuss their sub-optimality.

\subsection{Sub-optimality of $\chi^2$- and $G$-tests}
Practitioners frequently use $\chi^2$- and $G$-tests for conditional independence, which have nice asymptotic properties in fixed dimensional settings. In this subsection, we move beyond this fixed dimensional case and prove that these classical tests are markedly sub-optimal in terms of sample complexity. 
To define the $\chi^2$-test and $G$-test formally, let us write $o_{qrs} = \sum_{i=1}^n \mathds{1}(X_i = q) \mathds{1}(Y_i = r) \mathds{1}(Z_i = s)$ and $e_{qrs} = o_{q+s} o_{+rs} / o_{++s}$ where $o_{q+s} = \sum_{r \in [\ell_2]} o_{qrs}$, $o_{+rs} = \sum_{q \in [\ell_1]} o_{qrs}$ and $o_{++s} = \sum_{q \in [\ell_1],r \in [\ell_2]} o_{qrs}$, respectively, for $q \in [\ell_1], r \in [\ell_2], s \in [d]$. Given this notation, the $\chi^2$-test and $G$-test are based on the following test statistics
\begin{align} \label{Eq: chi and G statistics}
	\chi^2 = \sum_{q \in [\ell_1],r \in [\ell_2], s \in [d]} \frac{(o_{qrs} - e_{qrs})^2}{e_{qrs}} \quad \text{and} \quad G =  2 \sum_{q \in [\ell_1],r \in [\ell_2], s \in [d]} o_{qrs}  \log \biggl( \frac{o_{qrs} }{e_{qrs}} \biggr).
\end{align}
These test statistics converge to a $\chi^2$ distribution with $(\ell_1-1) \times (\ell_2-1)  \times d$ degrees of freedom under the null of conditional independence and under some regularity conditions \citep{tsamardinos2010permutation}. Based on this asymptotic result, $\chi^2$-test and $G$-test reject the null when $\chi^2$ and $G$ are larger than the $1-\alpha$ quantile of the $\chi^2$ distribution with $(\ell_1-1) \times (\ell_2-1)  \times d$ degrees of freedom. We first emphasize that these classical tests do not control the type I error uniformly over the null distributions and their validity guarantee requires that the sample size go to infinity. This is even true for the simplest case where $d=1$, which corresponds to the unconditional independence problem~\citep[see][for details]{berrett2021usp}. Moreover, their asymptotic power can be exactly equal to zero in some regimes where the sample size is much larger than the bound in \eqref{Eq: sample complexity} as shown below.

\begin{proposition}[Sub-optimality of $\chi^2$- and $G$-tests]  \label{Proposition: sub-optimality}
	Assume that $\ell_1 = \ell_2 = 2$, $\varepsilon = 0.25$ and $\alpha \in (0,1)$ is some fixed constant. Further assume that $d = n \times r_n$ where $r_n$ is an arbitrary positive sequence that increases to infinity as $n \rightarrow \infty$. In this scenario, the worse case power of $\chi^2$- and $G$-tests approach zero as 
	\begin{align*}
		\lim_{n \rightarrow \infty} \inf_{p \in \mathcal{P}_1(\varepsilon)}\mP_{p}(\chi^2 > q_{1-\alpha,d}) = 0 \quad \text{and} \quad \lim_{n \rightarrow \infty} \inf_{p \in \mathcal{P}_1(\varepsilon)}\mP_{p}(G > q_{1-\alpha,d}) = 0,
	\end{align*}
	where $q_{1-\alpha, d}$ denotes the $1-\alpha$ quantile of the $\chi^2$ distribution with $d$ degrees of freedom.
\end{proposition}

We provide some remarks on this result.

\begin{remark} \leavevmode \normalfont
	\begin{itemize}
		\item As shown in Theorem~\ref{Theorem: Multinomial Sampling using Permutation}, the proposed tests can achieve significant power (indeed rate optimal when $\ell_1=\ell_2=2$) under the same scenario and $r_n \lesssim n^{1/6}$. On the other hand,  $\chi^2$- and $G$-tests have asymptotically zero power for any $r_n$ that increases to infinity, which highlights sub-optimality of these classical tests. Moreover, the \texttt{UCI} and \texttt{wUCI} tests are valid over the entire class of null distributions unlike asymptotic $\chi^2$- and $G$-tests. 
		\item At a high-level, the reason behind this negative result is that the critical values of $\chi^2$- and $G$-tests are not adaptive to the underlying distribution. More specifically, we can think of a setting where most of conditional bins are empty with high probability, i.e.,~the intrinsic dimension of $Z$ is much smaller than $d$. In this case, it is more natural to use a critical value that reflects the intrinsic dimension rather than the ambient dimension. However, the $\chi^2$- and $G$-tests do not take this intuition into account, and their test statistics can be much smaller than $q_{1-\alpha,d}$ under the alternative. This leads to asymptotically zero power as we formally prove in Appendix~\ref{Proof: Proposition: sub-optimality}. 
		\item This issue can be alleviated by using the permutation approach where empty bins are ignored in calibration~\citep{tsamardinos2010permutation}. Nevertheless, it is unknown whether the permutation-based $\chi^2$- and $G$-tests are optimal in terms of sample complexity. We leave this important question for future work.
	\end{itemize}
\end{remark}

\section{Numerical analysis} \label{Section: Numerical analysis}
In this section, we provide numerical results that compare the proposed tests (\texttt{UCI}-test in Algorithm~\ref{Algorithm: Unweighted U-statistic} and \texttt{wUCI}-test in Algorithm~\ref{Algorithm: Weighted U-statistic without splitting}) with $\chi^2$-test and $G$-test under various scenarios. For a fair comparison, we calibrate both $\chi^2$-test and $G$-test using the permutation method (as in Algorithm~\ref{Algorithm: Unweighted U-statistic}) and reject the null when their permutation $p$-values are less than or equal to the significance level $\alpha$. Throughout our simulations, we set $\alpha = 0.05$ and the number of permutations $B = 199$. All the power values reported in this section are estimated by Monte Carlo simulation with 10,000 repetitions.

\subsection{Simulated data examples} \label{Section: Simulated data examples}

We start by comparing the power of the considered tests based on synthetic datasets. We only focus on the power results given that all of the tests are calibrated by the permutation procedure, resulting in valid type I error control in any finite sample sizes. There are eight different scenarios that we consider under the alternative where the domain sizes of $X,Y,Z$ are set to $\ell_1 = \ell_2 = 20$ and $d=10$, respectively. The considered scenarios are described as follows. 
\begin{itemize}
	\item \textbf{Scenario 1.} Set $p_Z$ to be uniform over $[d]$.  For each $z \in [d]$, (i) first let $p_{X,Y|Z}(x,y\,|\,z) \propto x^{-2} y^{-2}$, (ii) then replace $p_{X,Y|Z}(\ell_1,\ell_2\,|\,z)$ with $0.015$, and (iii) finally normalize $p_{X,Y|Z}$ to have its sum to be one. This setting results in a strong signal in $\chi^2$ divergence but relatively weaker signal in the $L_2$ distance over bins. 
	\item \textbf{Scenario 2.} Set $p_Z$ to be uniform over $[d]$. For each $z \in [d]$, (i) let $p_{X,Y|Z}(x,y\,|\,z) \propto x^{-2} y^{-2}$, (ii) set $\delta = \min\{p_{X,Y|Z}(x,y\,|\,z): x \in [2], y \in [2]\}$, and (iii) perturb $p_{X,Y|Z}$ by replacing $p_{X,Y|Z}(x,y\,|\,z)$ with $p_{X,Y|Z}(x,y\,|\,z) + (-1)^{x+y} \delta$ for $x \in [2]$ and $y \in [2]$. The resulting probability vector has a small signal in $\chi^2$ divergence, but relatively stronger signal in the $L_2$ distance over bins. 
	\item \textbf{Scenario 3.} Set $p_Z$ to be uniform over $[d]$. For each $z \in [d]$, (i) set $p_{X,Y|Z}(x,y\,|\,z) = 0$ for all $x \in [\ell_1]$ and $y \in [\ell_2]$, (ii) set $p_{X,Y|Z}(1,1\,|\,z) = (1-q)^2$, $p_{X,Y|Z}(1,y\,|\,z) = (1-q)q(\ell_1 - 1)^{-1}$ for $y \in [\ell_2] \! \setminus \! \{1\}$, $p_{X,Y|Z}(x,1\,|\,z) = (1-q)q (\ell_1 - 1)^{-1}$ for $x \in [\ell_1] \! \setminus \! \{1\}$, $p_{X,Y|Z}(x,y\,|\,z) = q^2 (\ell_1 - 1)^{-1}$ for $x = y \in [\ell_1] \! \setminus \! \{1\}$ where $q = 0.2$. This simulation setting is borrowed from \cite{zhang2022normal}.
	\item \textbf{Scenario 4.} Set $p_Z$ to be uniform over $[d]$. For each $z \in [d]$, let $p_{X,Y|Z}(x,y\,|\,z)  = \{1 + (-1)^{x+y}\} \ell_1^{-1} \ell_2^{-1}$ be a perturbed uniform distribution. This is the setting where $\chi^2$-test, $G$-test and \texttt{UCI}-test perform similarly for unconditional independence testing. See Figure 5 of \cite{berrett2021usp}. 
	\item \textbf{Scenario 5.} Set $p_Z$ to be uniform over $[d]$. For $z = 1$, let $p_{X,Y|Z}(x,y\,|\,z)  = 0.25$ for $x \in [2], y \in [2]$ and zero otherwise. In other words, there is no signal in the first category of $Z$. On the other hand, for $z \in [d] \! \setminus \! \{1\}$, set $p_{X,Y|Z}(x,y\,|\,z)  = \{1 + (-1)^{x+y}\} \ell_1^{-1} \ell_2^{-1}$ as in Scenario 4.
	\item \textbf{Scenario 6.} Set $p_Z$ to be uniform over $[d]$. For $z = 1$, $p_{X,Y|Z}(1,1\,|\,z) = p_{X,Y|Z}(2,2\,|\,z) = 0.4$, $p_{X,Y|Z}(1,2\,|\,z) =p_{X,Y|Z}(2,1\,|\,z) = 0.1$ and zero otherwise. On the other hand, for $z \in [d] \! \setminus \! \{1\}$, set $p_{X,Y|Z}(x,y\,|\,z)  = \ell_1^{-1} \ell_2^{-1}$, i.e.,~$X \, \indep \, Y$ for $z \in [d] \! \setminus \! \{1\}$, resulting in a sparse alternative.
	\item \textbf{Scenario 7.} Set $p_Z(z) \propto z^{-1}$. For $z \in [d]$, let $p_{X,Y|Z}(x,y\,|\,z)  = \{1 + (-1)^{x+y} z^{-1} \} \ell_1^{-1} \ell_2^{-1}$. By construction, the signal becomes weaker as $z$ increases and the sample size $\sigma_z$ tends to be smaller as $z$ increases. 
	\item \textbf{Scenario 8.} $p_Z(z) \propto z^{-1}$. For $z \in [d]$, let $p_{X,Y|Z}(x,y\,|\,z)  = \{1 + (-1)^{x+y} (d- z + 1)^{-1} \} \ell_1^{-1} \ell_2^{-1}$. Note that the signal becomes stronger as $z$ increases, and the sample size $\sigma_z$ tends to be smaller as $z$ increases. To put it in another way, this is the reverse setting of Scenario 7.
\end{itemize}

\begin{figure}[t!]
	\centering 
	\begin{tabular}{cc}
		\includegraphics[width=0.49\textwidth ]{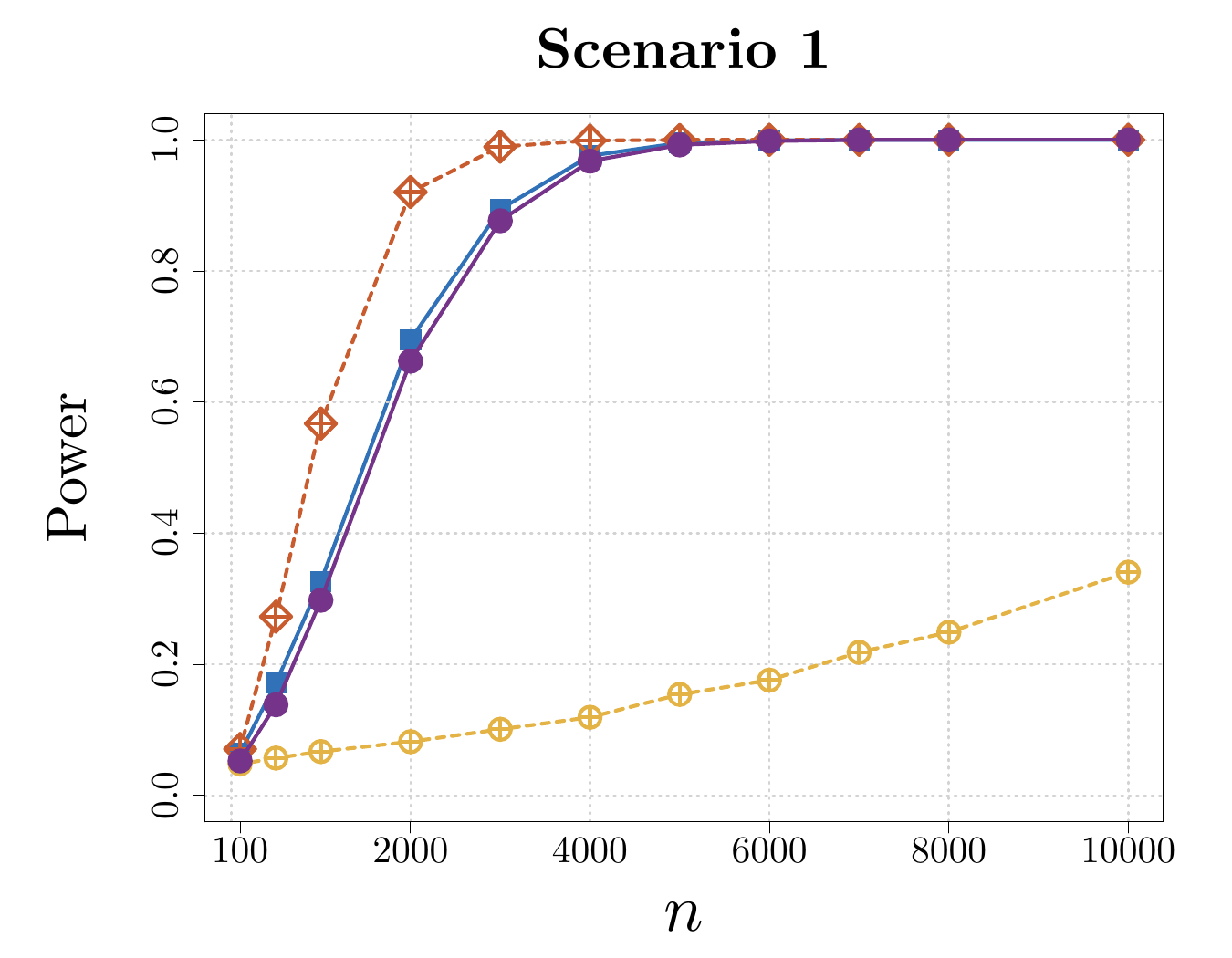} & 
		\includegraphics[width=0.49\textwidth ]{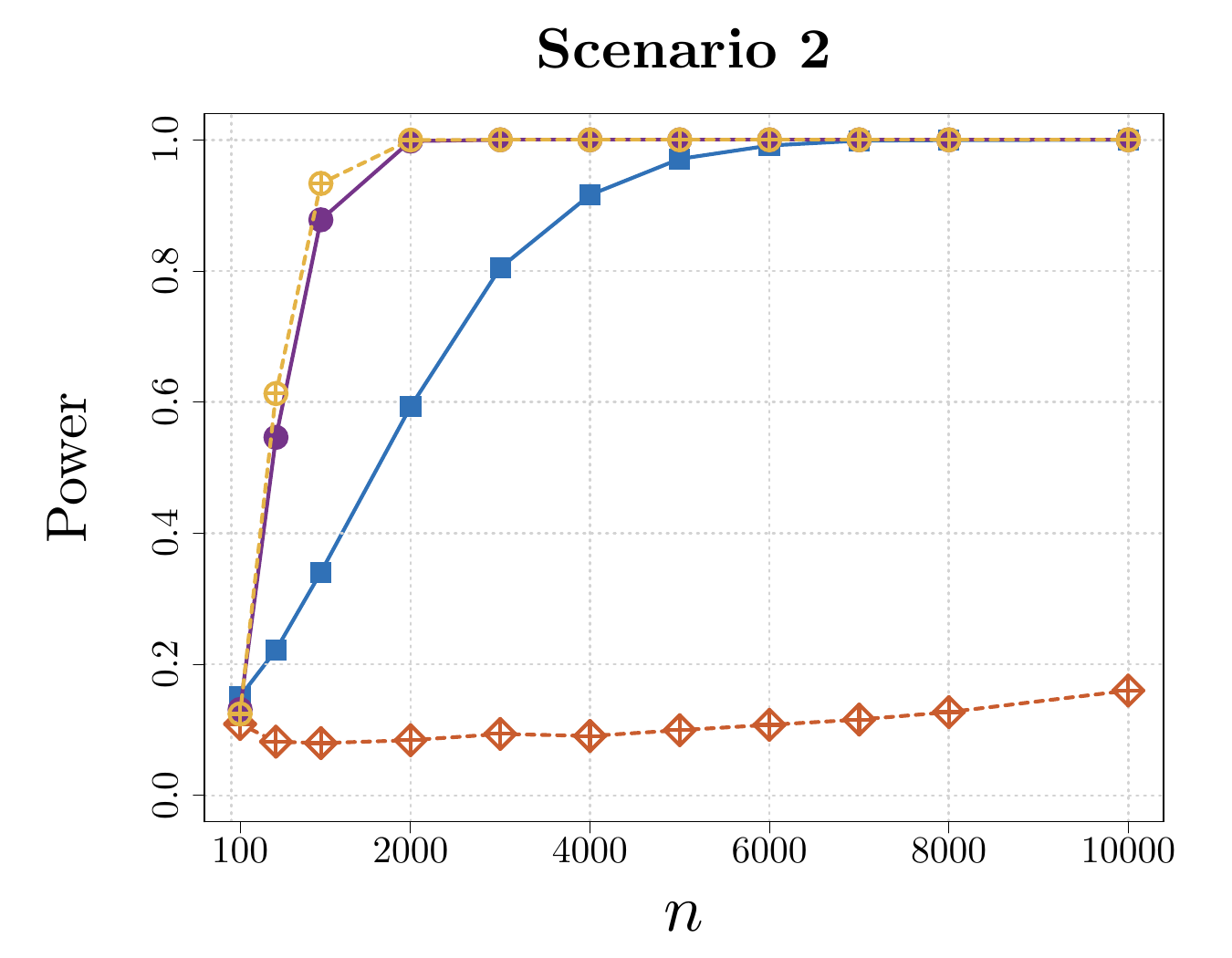} \\ 
		\includegraphics[width=0.49\textwidth ]{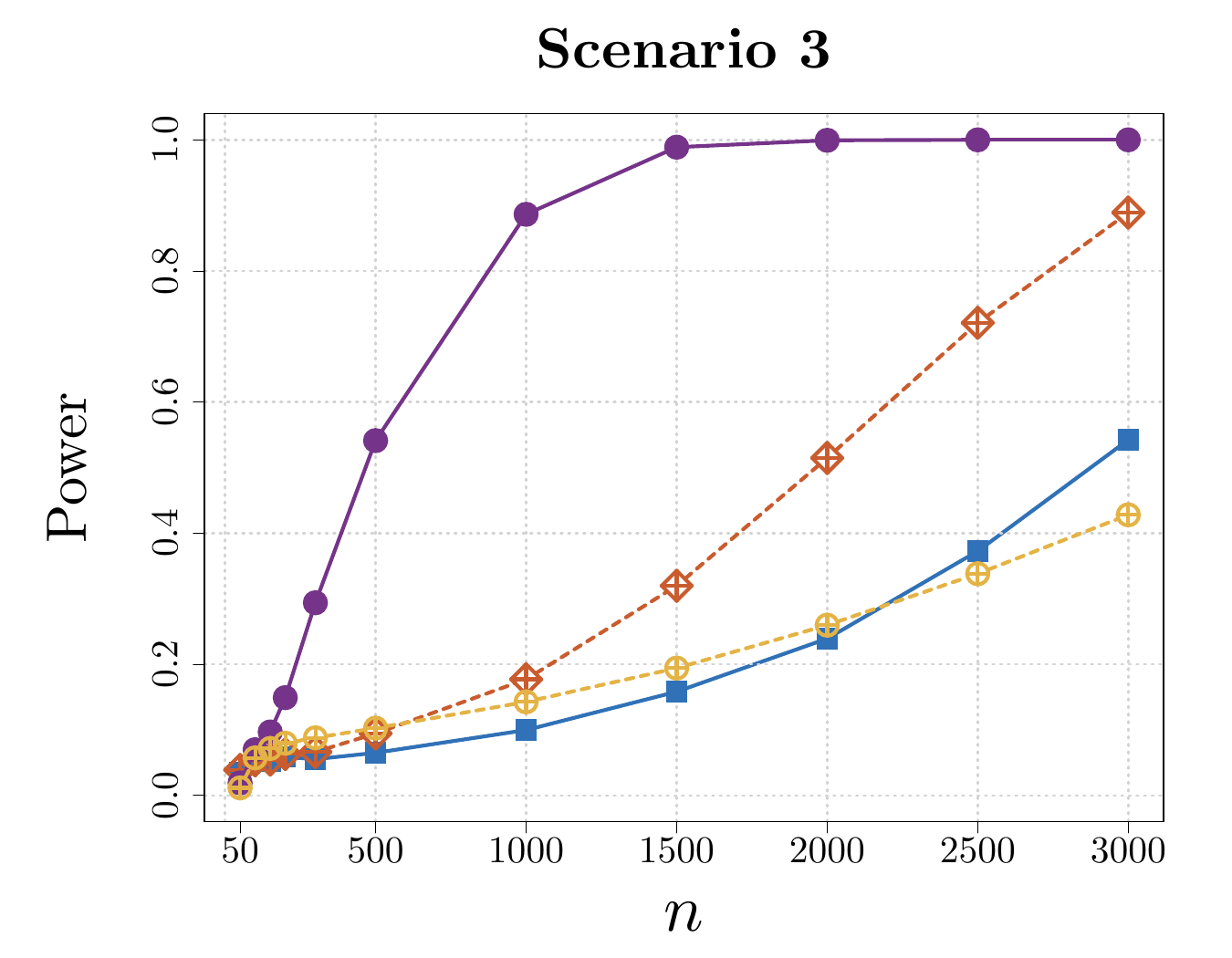} & 
		\includegraphics[width=0.49\textwidth ]{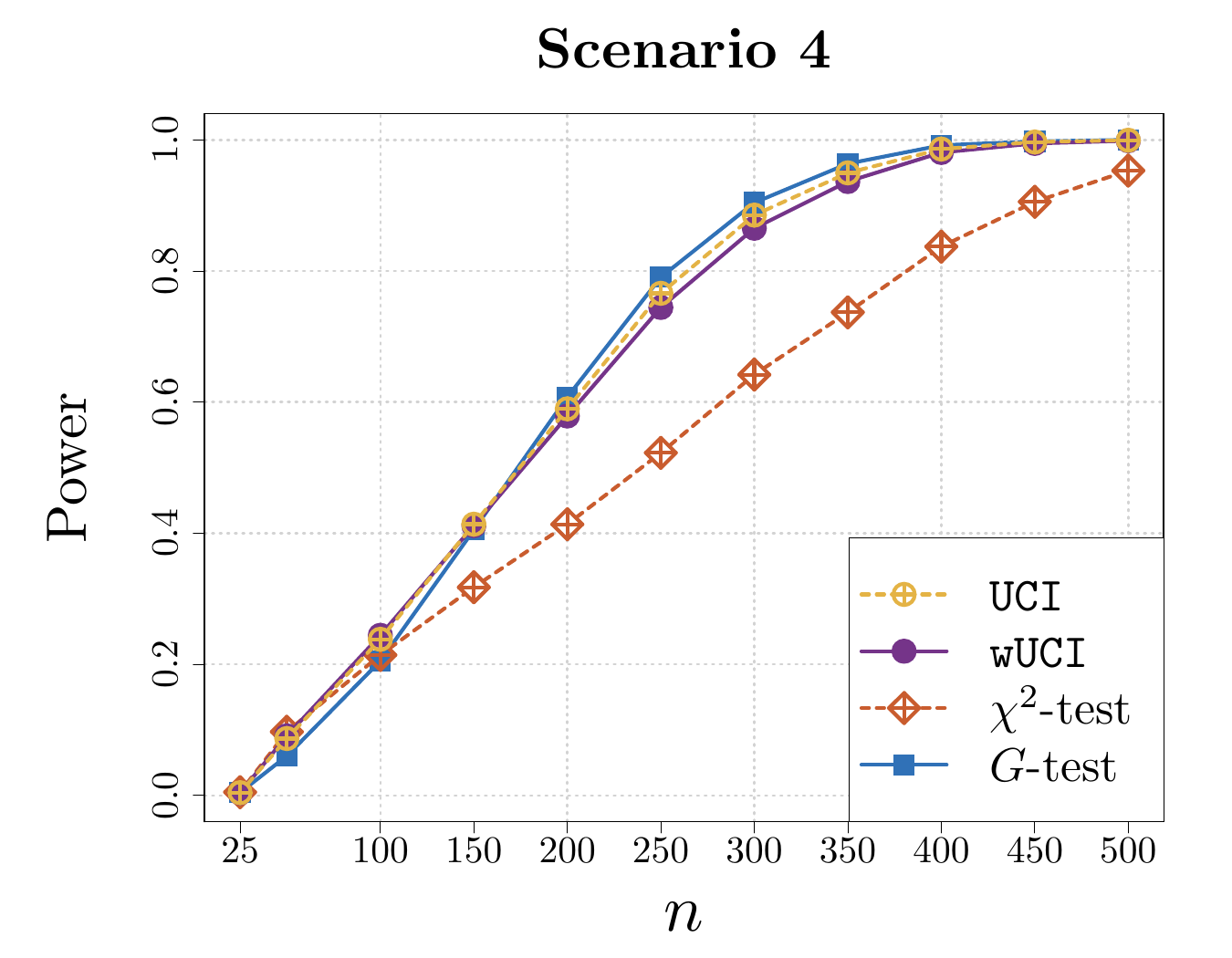}\\
	\end{tabular}
	\caption{Power comparisons of the considered tests in Scenario 1--Scenario 4 described in Section~\ref{Section: Simulated data examples}.} \label{Figure: Power I}
\end{figure}

\begin{figure}[t!]
	\centering 
	\begin{tabular}{cc}
		\includegraphics[width=0.49\textwidth ]{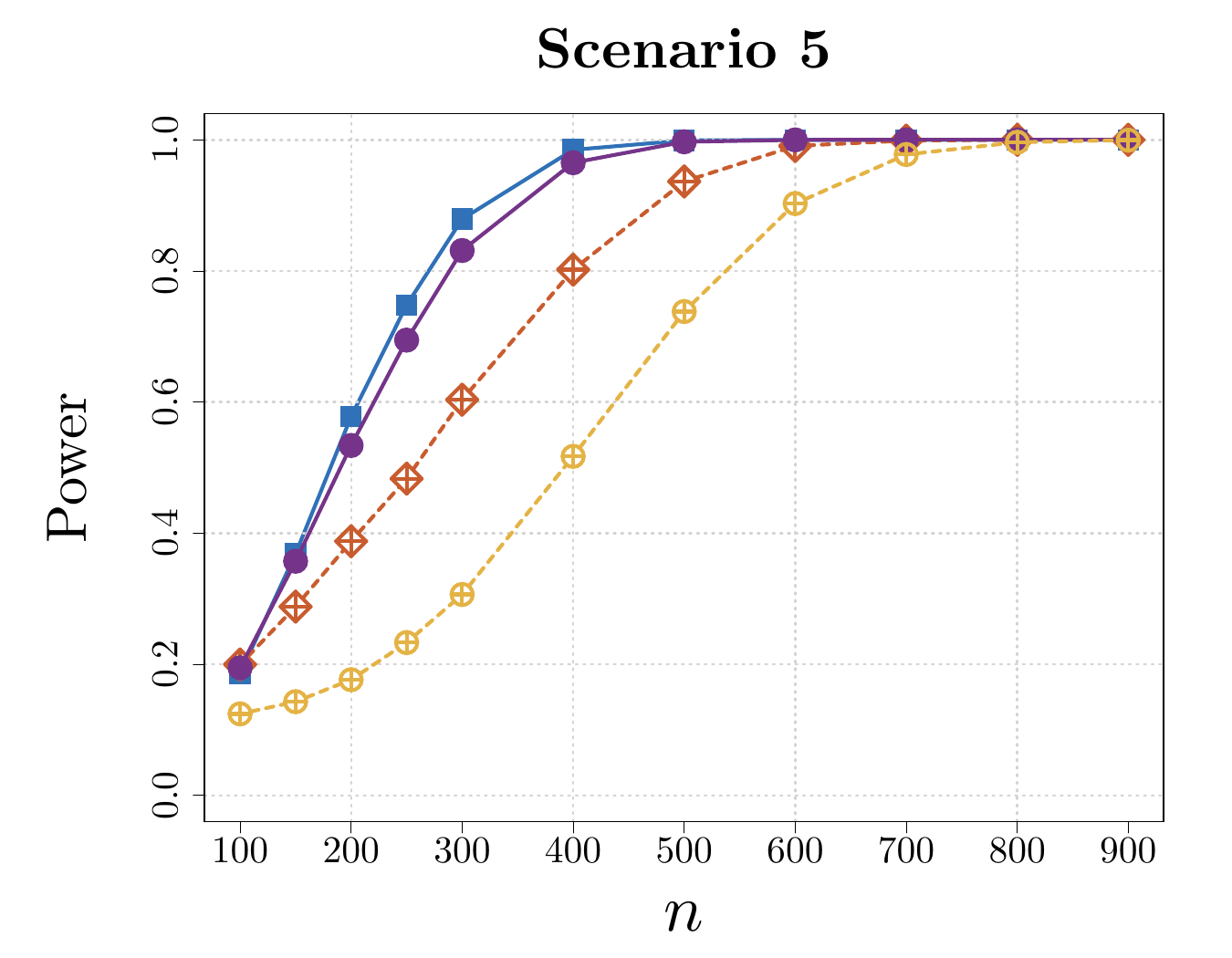} & 
		\includegraphics[width=0.49\textwidth ]{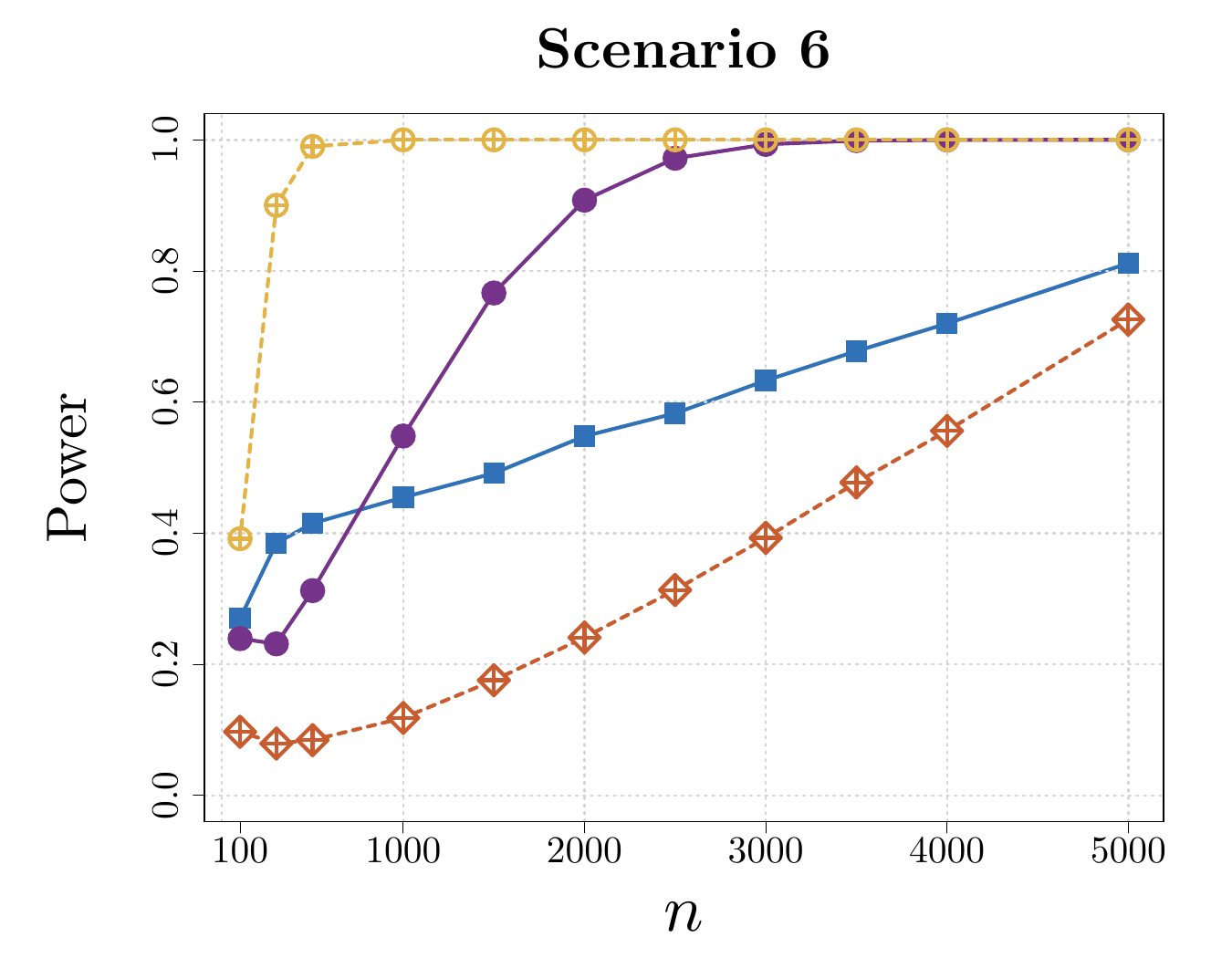} \\ 
		\includegraphics[width=0.49\textwidth ]{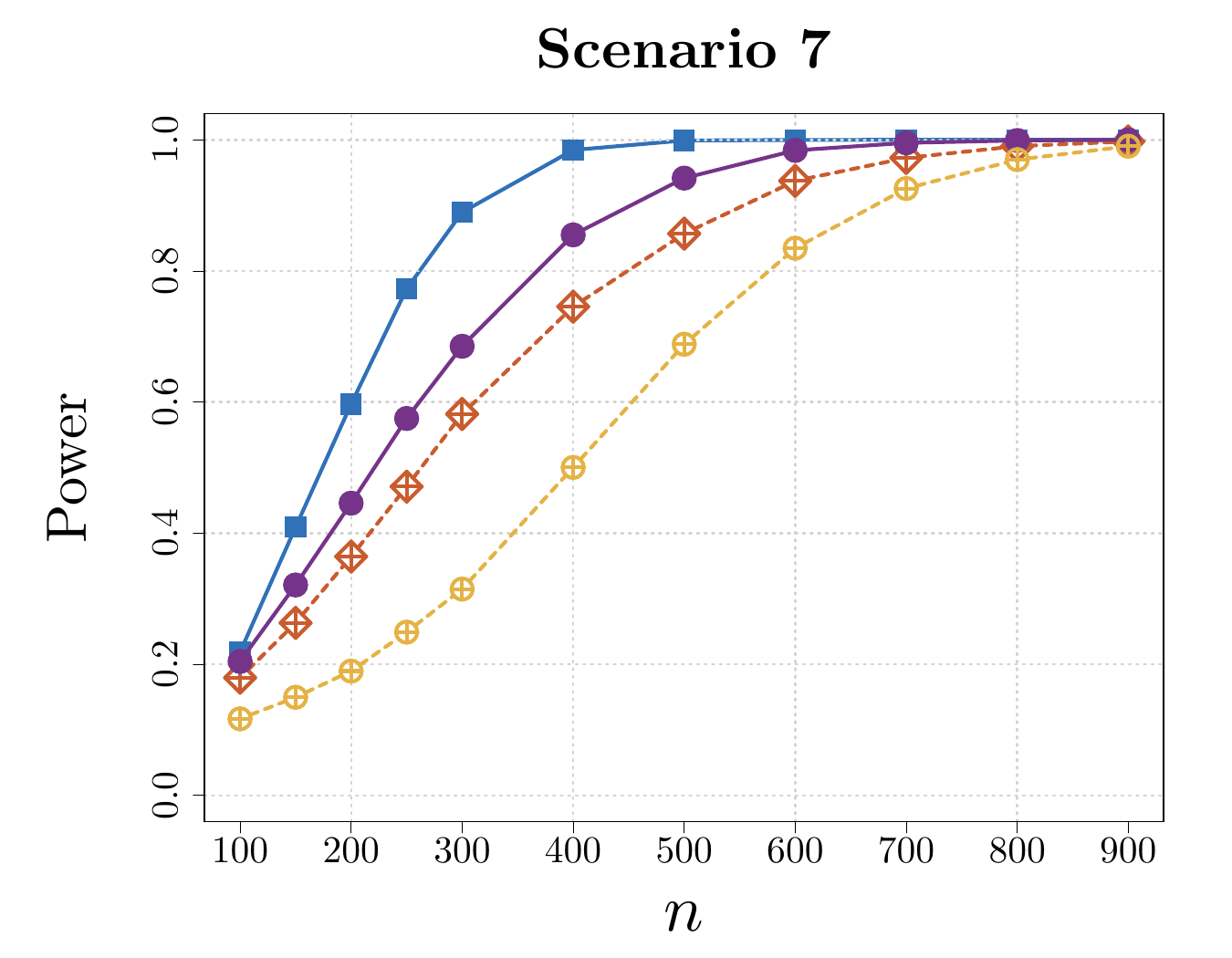} & 
		\includegraphics[width=0.49\textwidth ]{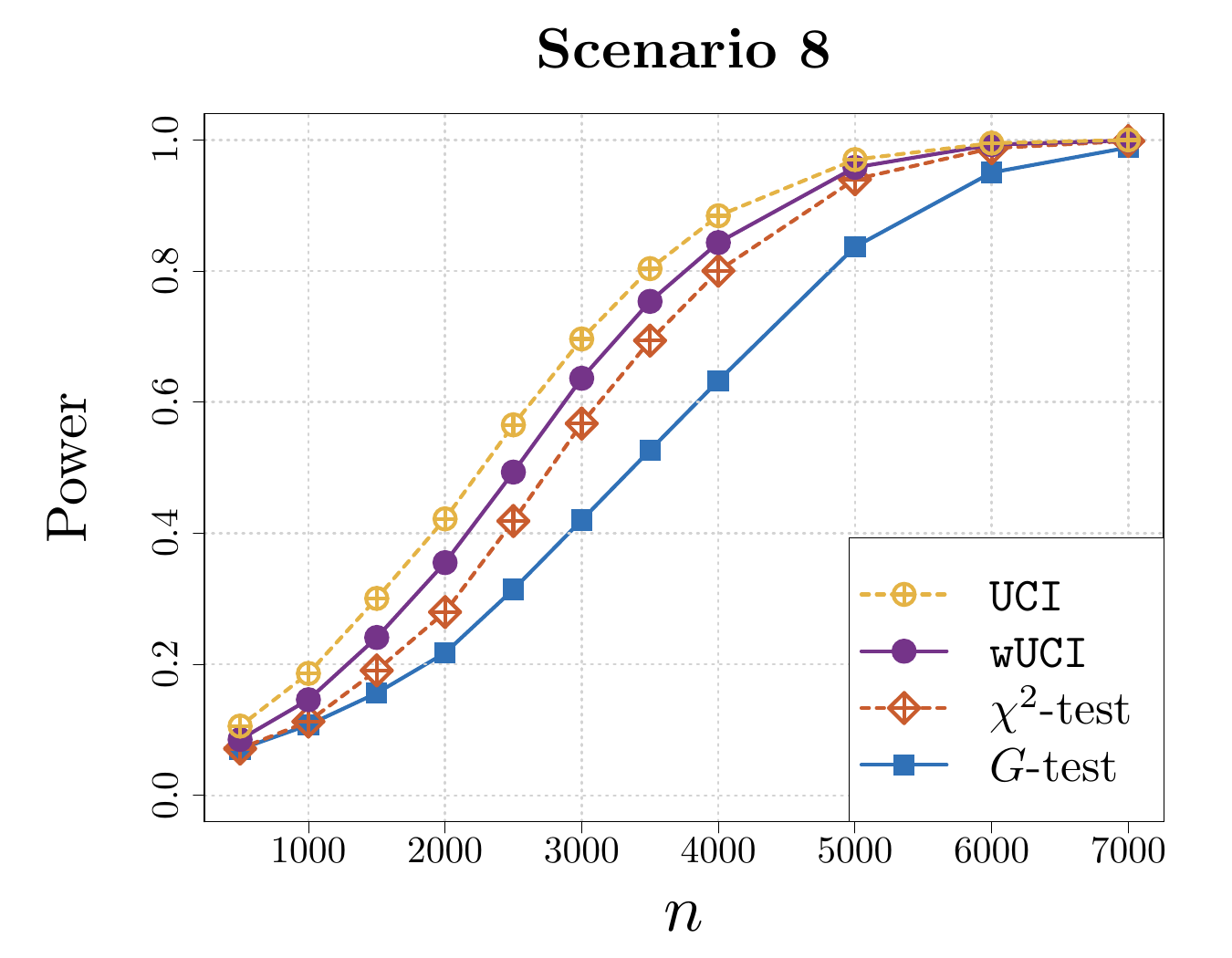}\\
	\end{tabular}
	\caption{Power comparisons of the considered tests in Scenario 5--Scenario 8 described in Section~\ref{Section: Simulated data examples}.} \label{Figure: Power II}
\end{figure}

The results are presented in Figure~\ref{Figure: Power I} and Figure~\ref{Figure: Power II}. It is clear from the results that no test outperforms the others over all scenarios. In particular, $\chi^2$-test has the highest power when the underlying distribution has a strong signal in $\chi^2$ divergence such as Scenario 1, and similarly, \texttt{UCI}-test performs well when there is a strong signal in the $L_2$ distance such as Scenario 2. It is also interesting to observe that \texttt{wUCI}-test shows an impressive performance compared to the others in Scenario 3, and all of the tests perform similarly in Scenario 4 (except $\chi^2$-test with large sample sizes) where the corresponding null distribution of $(X,Y) \,|\,Z$ is uniform over bins. In Scenario 5 and Scenario 6, we are essentially in a situation where $\ell_1=\ell_2=2$ for $z=1$ and $\ell_1 = \ell_2 = 20$ for $z \in [d] \! \setminus \! \{1\}$. In these scenarios, the first bin of $Z$ plays a less important role than the other bins in determining the value of $\chi^2$ and $G$ statistics (these statistics have higher variance when $\ell_1$ and $\ell_2$ are large). In contrast, the test statistic for the \texttt{UCI}-test is mostly dominated by the data from the first bin of $Z$ in the same scenarios. This explains the outstanding performance of $\chi^2$-test and $G$-test in Scenario 5 where signals are spread out over bins except the first one. On the other hand, $\chi^2$-test and $G$-test have low power in Scenario 6 where only the first bin of $Z$ has a signal. Under the same scenarios, \texttt{UCI}-test behaves in the opposite way, attaining high power when the first bin is significant. Another interesting observation is that \texttt{wUCI}-test performs as powerful as $G$-test in Scenario 5 whereas it outperforms both $\chi^2$- and $G$-tests in Scenario 6 when the sample size is large. This may be explained by the choice of weights in its statistic that roughly interpolate $\chi^2$ weights and uniform weights as explained in Remark~\ref{Remark: discussion on weights}. We also note that the test statistic for \texttt{UCI}-test is a linear combination of U-statistics weighted by sample sizes over bins. This explains the relatively lower power of \texttt{UCI}-test than $\chi^2$- and $G$-tests in Scenario 7 where the bins with smaller sample sizes tend to have stronger signals. In contrast, we observe the opposite behavior in Scenario 8. On the other hand, \texttt{wUCI} performs the second best in both Scenario 7 and Scenario 8.

To summarize, we observe that different tests perform better than the others under different scenarios. The proposed tests often dominate the classical ones when there are strong signals, especially in the $L_2$ distance, over bins with large sample sizes. On the other hand, it is possible to design situations such as Scenario 1 where the proposed tests attain lower power than the classical ones. Nevertheless, the classical tests, especially $\chi^2$-test, can fail badly in terms of the worse-case performance.  By contrast, our proposals demonstrate robust performance across different scenarios, indicating that they can work as practical tools that complement classical $\chi^2$-test and $G$-test under various scenarios.

\subsection{Real-world data examples}
We next provide numerical illustrations based on real-world datasets. 

\vskip .5em

\noindent \textbf{Admission dataset.} The first dataset that we look at is the Berkeley admissions dataset, which is a well-known example of Simpson's paradox~\citep{bickel1975sex}. As summarized in Table~\ref{Table: Admissons}, the dataset consists of 4,526 applications with 3 variables $(X,Y,Z)$ where $X$ and $Y$ are binary variables, representing the gender (male or female) and the admission status (admitted or rejected), respectively. The conditional variable $Z$ takes the department name among $\{$A, B, C, D, E, F$\}$. When the dataset is aggregated over the departments, it appears that male applicants are more likely to be admitted than woman applicants. However, as reported by \cite{bickel1975sex}, there seems to exist a bias in favor of women when looking at the individual departments, indicating the existence of conditional dependence. We assess this claim of conditional dependence by implementing the considered permutation tests. For this dataset, the corresponding $p$-values are computed as 0.005 for both $\chi^2$- and $G$-tests, 0.04 for \texttt{UCI}-test and 0.03 for \texttt{wUCI}-test, respectively. All the $p$-values are significant at level $\alpha = 0.05$, revealing evidence of conditional dependence.

\begin{table}[h!]
	\begin{center} 	\caption{Admissions data at University of California, Berkeley from the six largest departments in 1973} \label{Table: Admissons}
	\setlength{\tabcolsep}{10pt}
	\renewcommand{\arraystretch}{1.1}
	\begin{tabular}{ccccc}
		\toprule
		\multicolumn{1}{l}{} & \multicolumn{2}{c}{Men} & \multicolumn{2}{c}{Women} \\  \cmidrule{2-3} \cmidrule{4-5}
		Major              & Applicants  & Admitted  & Applicants   & Admitted   \\   \midrule
		A                    & 825         & 62        & 108          & 82         \\ 
		B                    & 560         & 63        & 25           & 68         \\ 
		C                    & 325         & 37        & 593          & 34         \\ 
		D                    & 417         & 33        & 375          & 35         \\ 
		E                    & 191         & 28        & 393          & 24         \\ 
		F                    & 373         & 6         & 341          & 7          \\ \bottomrule
	\end{tabular}
	\end{center}
\end{table}

\vskip 1em

\begin{figure}[h!]
	\centering 
	\begin{tabular}{cc}
		\includegraphics[width=0.49\textwidth ]{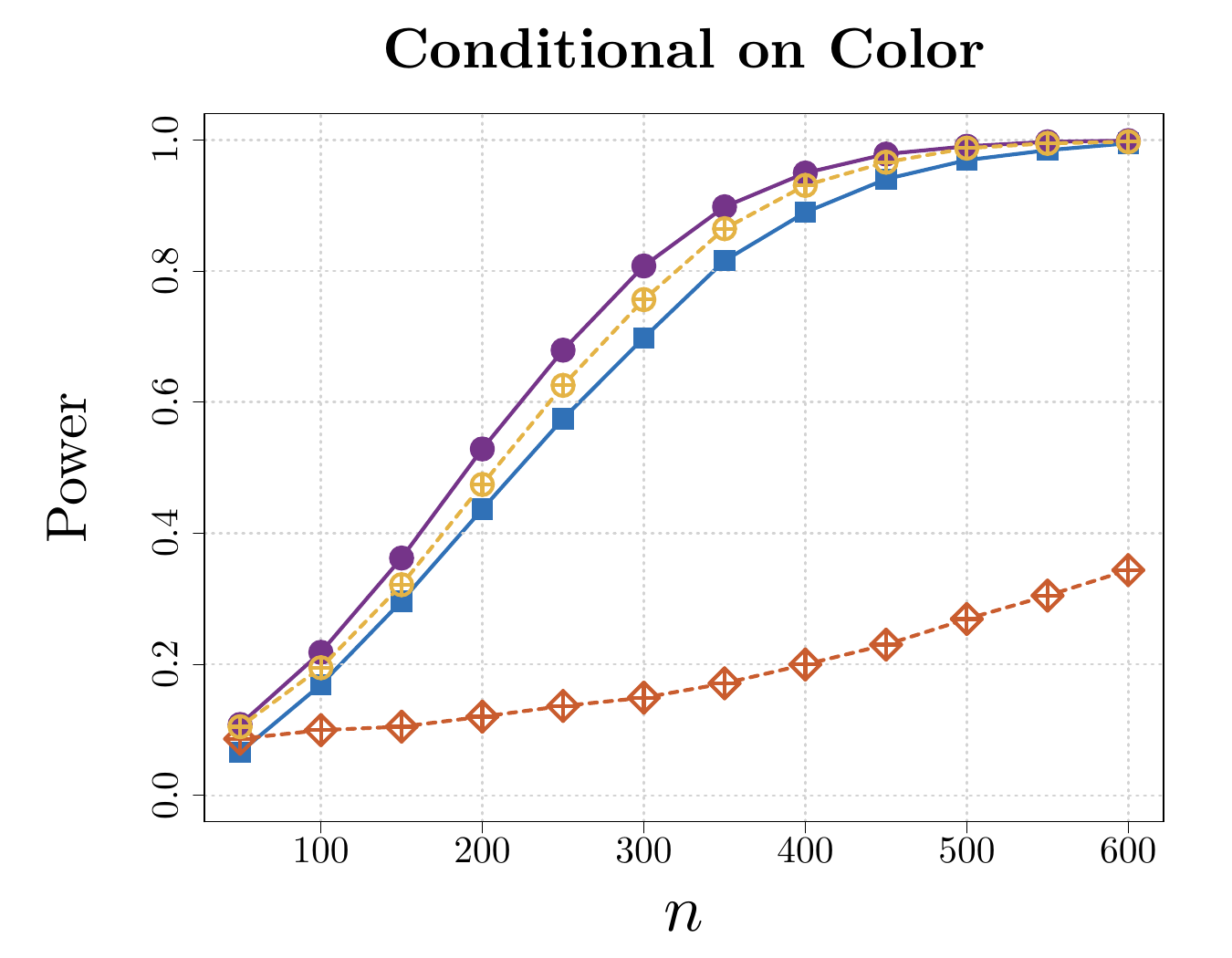} & 
		\includegraphics[width=0.49\textwidth ]{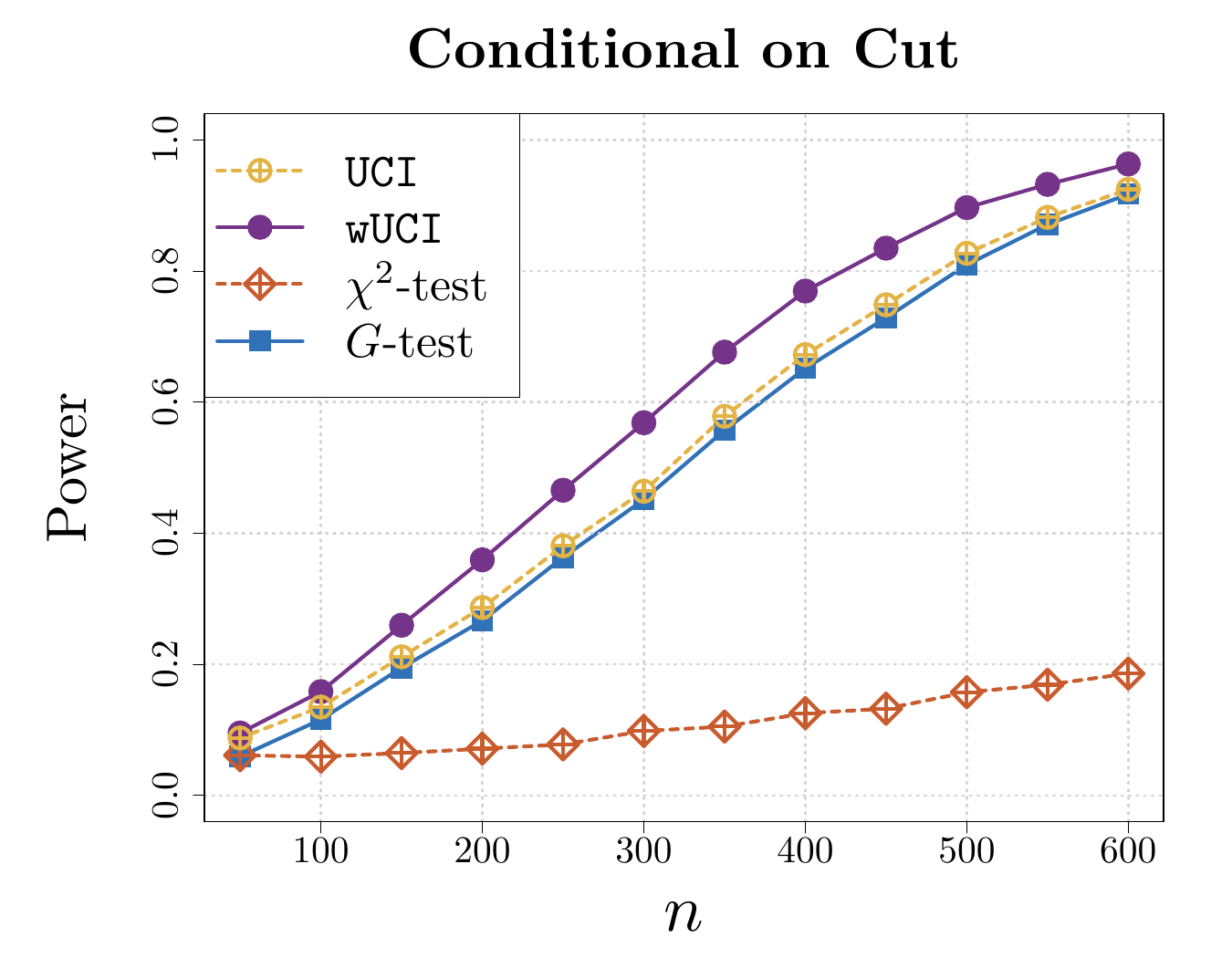} \\
	\end{tabular}
	\caption{Power comparisons of the considered tests based on the diamonds dataset. Both panels analyze independence between the (categorized) price and clarity variables conditional on the color variable and the cut variable, respectively. All of the tests have increasing power as the sample size increases. Markedly, $\chi^2$-test has significantly lower power than the other tests, whereas \texttt{wUCI}-test seems to perform the best for this dataset.} \label{Figure: Diamonds}
\end{figure}

\noindent \textbf{Diamonds dataset.} Next we consider the diamonds dataset available in \texttt{R} package \texttt{ggplot2}. The dataset contains the information of 53,940 diamonds including their price, clarity, color, quality of the cut, etc. In our analysis, the price variable is partitioned into 100 intervals of equal size. We set the corresponding categorized price variable as $X$ and set the clarity variable as $Y$. The clarity variable has 8 categories (I1, SI2, SI1, VS2, VS1, VVS2, VVS1, IF) and it measures the purity of a diamond. The conditional variable $Z$ is chosen to be either the cut variable or the color variable in our analysis. Both variables are discrete with 5 (Fair, Good, Very Good, Premium, Ideal) and  7 (D, E, F, G, H, I, J) categories, respectively. In the experiments, we treat the entire dataset as the population (thereby the ground truth is known to us)\footnote{We verified numerically that $X$ and $Y$ are conditionally dependent on $Z$ by comparing the distribution of $p_{X,Y,Z}$ and $p_{X|Z}p_{Y|Z}p_Z$.} from which we randomly draw $n$ observations without replacement. Based on this subsample of size $n$, we compute the permutation $p$-values of the considered tests and we repeat this process $10,000$ times to estimate their power. The results can be found in Figure~\ref{Figure: Diamonds} where we collect the power of $\{\texttt{UCI}, \texttt{wUCI}, \chi^2, G\}$-tests by changing the sample size $n$. The left panel of Figure~\ref{Figure: Diamonds} provides the power results when the conditional variable is set to be the color variable. As can be seen, $\chi^2$-test has the significantly lower power than the others. Among the other three tests, \texttt{wUCI}-test has the highest power followed by \texttt{UCI}-test while the difference is minor. We can see a similar pattern from the right panel of Figure~\ref{Figure: Diamonds} where the conditional variable is set to be the cut variable. These results highlight the practical value of the proposed tests in analyzing real-world datasets where classical tests potentially suffer from low power.

\section{Discussion} \label{Section: Discussion}
In this paper, we have revisited recent developments of CI testing for discrete data. Despite attractive theoretical properties, these recent tests have limited practical value, relying on Poissonization and unspecified constants in their critical values. In this work, we have made an attempt to bridge the gap between theory and practice by removing Poissonization and utilizing the Monte Carlo permutation method to calibrate test statistics. We have also complemented our theoretical results with a thorough numerical analysis and demonstrated certain benefits of the proposed tests over classical $\chi^2$- and $G$-tests. Finally, we have developed \texttt{R} package \texttt{UCI} that implements the proposed methods.

Our work leaves several important avenues for future research. One prominent direction is to depoissonize other sample complexity results in the literature using the tools developed in this paper. For instance, one can reproduce the results of \cite{neykov2021minimax,kim2021local} for continuous CI testing without Poissonization. Another direction which may be fruitful to pursue is to devise a CI test that incorporates prior information about potential alternative distributions. For example, suppose that we are in an alternative setting where only a handful of conditional categories are significant. In this case, it is possible to obtain a substantial power gain by using sparse weights in the proposed statistics, and one could analyze the resulting tests. Additionally, it would be interesting to see whether \texttt{wUCI}-test or other tests without sample splitting can achieve the sample complexity~\eqref{Eq: general sample complexity} in general. We leave these interesting questions for future work. 

\paragraph{Acknowledgements.} 
This work was partially supported by funding from the NSF grants DMS-2113684 and DMS-2310632, as well as an Amazon AI and a Google Research Scholar Award to SB. MN acknowledges support from the NSF grant DMS-2113684. IK acknowledges support from the Basic Science Research Program through the National Research Foundation of Korea (NRF) funded by the Ministry of Education (2022R1A4A1033384), and the Korea government (MSIT) RS-2023-00211073.

\bibliographystyle{apalike}
\bibliography{reference}

\clearpage 

\appendix

\section{Algorithm} \label{Section: Algorithm}
The below is the algorithm for \texttt{wUCI}-test used in our simulation study. 
\begin{algorithm} \caption{\texttt{wUCI}: weighted U-statistic permutation CI test} \label{Algorithm: Weighted U-statistic without splitting}
	\textbf{Input:} Sample~$\{(X_i,Y_i,Z_i)\}_{i=1}^n$, the number of permutations~$B$, significance level $\alpha$
	\vskip .5em
	\textbf{For} $j \in [B]$ \textbf{do}
	\vskip .5em
	\hskip 2em \textbf{For} $m \in [d]$ \textbf{do}
	\begin{itemize}[leftmargin=1.8cm]
		\item Generate $\pi \sim \mathrm{Uniform}(\Pi_{\sigma_m})$ independent of everything else.
		\item Compute $U_W^{\boldsymbol{b}}(D_m^\pi)$ as in (\ref{Eq: weighted U without splitting}) based on the permuted dataset $D_m^\pi$. 
	\end{itemize}
	\hskip 2em \textbf{End}
	\vskip .5em
	\hskip 2em Set $T_{W,j}^\dagger \leftarrow \sum_{m \in [d]} \mathds{1}(\sigma_m \geq 4) \sigma_m \omega_m U_W^{\boldsymbol{b}}(D_m^\pi)$.
	\vskip .5em
	\textbf{End}
	\begin{itemize}
		\item Set $T_W^\dagger \leftarrow \sum_{m \in [d]} \mathds{1}(\sigma_m \geq 4) \sigma_m \omega_m U_W^{\boldsymbol{b}}(D_m)$ computed without permutations.
		\item Compute the permutation $p$-value
		\begin{align*}
			p_{\mathrm{perm}} = \frac{1}{B+1} \Biggl[ \sum_{j=1}^B \mathds{1}(T_{W,j}^\dagger \geq T_W^\dagger) + 1 \Biggr].
		\end{align*}
	\end{itemize}
	\vskip .5em
	\textbf{Output:} Reject $H_0$ if $p_{\mathrm{perm}} \leq \alpha$; otherwise, accept $H_0$. 
\end{algorithm}

\section{Lemmas} \label{Section: Lemmas}
In this section, we prove several lemmas in \cite{canonne2018testing} under multinomial sampling. These results are main building blocks that lead to Theorem~\ref{Theorem: Multinomial Sampling} and Theorem~\ref{Theorem: Multinomial Sampling using Permutation} in the main text. We start by showing that the result in Lemma 3.1 of \cite{canonne2018testing} holds without Poissonization.

\begin{lemma} \label{Lemma: lower bound}
	Suppose that $X \sim \mathrm{Binomial}(n,p)$ where $n \geq 4$. Then there exists an absolute constant $\gamma >0.85$ such that 
	\begin{align*}
		\mE[X \mathds{1}(X \geq 4)] \geq \gamma \min\{np, (np)^4\}.
	\end{align*}
\end{lemma}

The second lemma that we prove corresponds to Claim 2.1 of \cite{canonne2018testing}.

\begin{lemma} \label{Lemma: Variance bound}
	Suppose that $X \sim \mathrm{Binomial}(n,p)$ where $n \geq 4$. Then there exists a constant $\gamma >0$ such that 
	\begin{align*}
		\mathrm{Var}[X \mathds{1}(X \geq 4)] \leq \gamma \mE[X \mathds{1}(X \geq 4)].
	\end{align*}
\end{lemma}

The following lemma proves Claim 2.2 of \cite{canonne2018testing} for a binomial random variable with parameters $n,p$. It is worth pointing out that the the original statement of Claim 2.2 in \cite{canonne2018testing} contains an error, which has been corrected by \cite{kim2022comments}.

\begin{lemma} \label{Lemma: Variance bound Claim 22}
	Suppose that $X \sim \mathrm{Binomial}(n,p)$ where $n \geq 9$ and let $a,b \geq 2$. Then there exists a constant $C>0$ such that
	\begin{align*}
		&\mathrm{Var}\bigl[ X\sqrt{\min\{X,a\}\min\{X,b\}} \mathds{1}(X\geq 4) \bigr] \\[.5em] 
		\leq & ~ C \min\{np,\sqrt{ab}\} \mE\bigl[ X \sqrt{\min\{X,a\}\min\{X,b\}} \mathds{1}(X\geq 4) \bigr].
	\end{align*}
\end{lemma}

The following lemma proves Claim 2.3 of \cite{canonne2018testing} for a binomial random variable with parameters $n,p$.

\begin{lemma} \label{Lemma: expectation lower bound}
	Suppose that $X \sim \mathrm{Binomial}(n,p)$ where $n \geq 9$. Then there exists a constant $C>0$ such that for any $a,b \geq 2$ and $\lambda = np$,
	\begin{align*}
		\mE\bigl[X \sqrt{\min\{X,a\}\min\{X,b\}} \mathds{1}(X\geq 4) \bigr] \geq C \min\bigl\{ \lambda\sqrt{\min(\lambda,a)\min(\lambda,b)}, \lambda^4 \bigr\}.
	\end{align*}
\end{lemma}

We say that random variables $X_1,\ldots,X_k$ are negatively associated if 
\begin{align*}
	\mathrm{Cov}\{f_1(X_i, i \in A_1), f_2(X_j, j \in A_2)  \} \leq 0, 
\end{align*}
for all increasing functions $f_1,f_2$ and for every pair of disjoint subsets $A_1,A_2$ of $\{1,2,\ldots,k\}$. It is well-known that a multinomial sample is negatively associated, which we record below.

\begin{lemma}[\cite{joag1983negative}] \label{Lemma: NA of multinomials}
	 Let $(\sigma_1,\ldots,\sigma_d)$ have a multinomial distribution with parameters $n$ and $(p_1,\ldots,p_d)$. Then for every pair of disjoint subsets $A_1,A_2$ of $\{1,\ldots,d\}$, it holds that
	 \begin{align*}
	 	\mathrm{Cov}\{f_1(\sigma_i, i \in A_1), f_2(\sigma_j, j \in A_2)\} \leq 0,
	 \end{align*} 
	 whenever $f_1$ and $f_2$ are increasing. 
\end{lemma}
The next lemma builds on the negative association of a multinomial random vector, which is proved in Appendix~\ref{Section: Proof of Lemma: NA}. 
\begin{lemma} \label{Lemma: NA}
	Suppose that $(\sigma_1,\ldots,\sigma_d)$ has a multinomial distribution with parameters $n$ and $(p_1,\ldots,p_d)$. Then for any $1 \leq i \neq j \leq d$
	\begin{align*}
		\mathrm{Cov}\bigl\{\sigma_i \mathds{1} (\sigma_i \geq 4), \sigma_j \mathds{1} (\sigma_j \geq 4)\bigr\} \leq 0.
	\end{align*}
	Moreover, let $\omega_i = \sqrt{ \min\{\sigma_i,\ell_1\} \min\{\sigma_i,\ell_2\}}$ for $\ell_1, \ell_2 \geq 0$. Then for any $1 \leq i \neq j \leq d$
	\begin{align*}
		\mathrm{Cov}\bigl\{\sigma_i \mathds{1} (\sigma_i \geq 4) \omega_i, \sigma_j \mathds{1} (\sigma_j \geq 4)\omega_j \bigr\} \leq 0.
	\end{align*}	
\end{lemma}

As shown in Lemma 6 of \cite{kim2021local}, rejecting the null when the permutation $p$-value is less than or equal to $\alpha$ is equivalent to rejecting the null when the corresponding test statistic is greater than the $1-\alpha$ quantile of the permutation distribution. We state this result below for completeness, and its proof can be found in Appendix~\ref{Section: Proof of Lemma: quantile}. 

\begin{lemma}[Quantile] \label{Lemma: quantile}
	Let $q_{1-\alpha}$ be the $1-\alpha$ quantile of the empirical distribution of $V,V_1,\ldots,V_B$. Then it holds that 
	\begin{align*}
		\mathds{1}\Biggl(  \frac{1}{B+1} \Biggl[ \sum_{i=1}^B \mathds{1}(V_i \geq V) +1 \Biggr] > \alpha \Biggr)  = \mathds{1}(V \leq q_{1-\alpha}). 
	\end{align*}
\end{lemma}

Finally, we state the Efron--Stein inequality for completeness \citep[e.g.,~Theorem 3.1 of][]{boucheron2013concentration}.
\begin{lemma}[Efron--Stein inequality] \label{Lemma: Efron-Stein}
	Let $X_1,\ldots, X_n$ be independent random variables and let $Z = f(\boldsymbol{X})$ be a square-integrable function of $\boldsymbol{X} = (X_1,\ldots,X_n)$. Moreover, $X_1',\ldots,X_n'$ are independent copies of $X_1,\ldots,X_n$ and define 
	\begin{align*}
		Z'_i = f(X_1,\ldots,X_{i-1},X_i',X_{i+1},\ldots,X_n).
	\end{align*} 
	Then 
	\begin{align*}
		\mathrm{Var}[Z] \leq \frac{1}{2} \sum_{i=1}^n \mE \bigl[ (Z - Z_i')^2 \bigr].
	\end{align*}
\end{lemma}

\section{Proofs} \label{Section: Proofs}
This section collects the proof of the main theorems and lemmas. 

\subsection{Proof of Theorem~\ref{Theorem: Multinomial Sampling}}
Given the lemmas in Appendix~\ref{Section: Lemmas} together with the results in \cite{canonne2018testing}, the claims in Theorem~\ref{Theorem: Multinomial Sampling} are almost immediate. We present more details below. Since the permutation test controls the type I error in any finite sample sizes, we only focus on the alternative case where the underlying distribution is $\varepsilon$-far from the null in the $L_1$ distance, i.e.,~$\inf_{q \in \mathcal{P}_0} \|p - q\|_1 \geq \varepsilon$, for both $T$ and $T_W$. 

\vskip 1em

\subsubsection{Claim for $T$} 

Starting with the first claim, we analyze the expectation and the variance of $T$ under multinomial sampling. Conditional on $\sigma_m$, $U(D_m)$ is an unbiased estimator of $\delta_m^2:= \|p_{X,Y\,|\,Z=m} - p_{X|Z=m}p_{Y|Z=m}\|_2^2$. Moreover, since $\sigma_m \sim \mathrm{Binomial}\bigl(n,p_Z(m)\bigr)$, Lemma~\ref{Lemma: lower bound} proves that
\begin{align*}
	\mE[T] ~=~ &  \sum_{m \in [d]} \mE\bigl[\mathds{1}(\sigma_m \geq 4) \sigma_m U(D_m) \bigr] \\[.5em]
	 = ~ &  \sum_{m \in [d]} \mE[\mathds{1}(\sigma_m \geq 4) \sigma_m] \cdot  \delta_m^2 \\[.5em]
	 \gtrsim ~ &  \sum_{m \in [d]}  \min\bigl\{np_Z(m), \bigl(np_Z(m)\bigr)^4\bigr\} \cdot \delta_m^2.  
\end{align*}	
Set the sample size $n = \beta' \max\bigl( \sqrt{d}/\varepsilon^{\prime 2}, \min (d^{7/8}/\varepsilon^\prime, d^{6/7}\varepsilon^{\prime 8/7}) \bigr)$ where $\varepsilon^\prime := \varepsilon/\sqrt{\ell_1 \ell_2}$ and $\beta' \geq 1$ is a sufficiently large constant. Using the inequality, we follow the same lines of the proof of Proposition 3.1 in \cite{canonne2018testing} and show that 
\begin{align*}
	\mE[T] ~ \gtrsim ~ \beta' \sqrt{\min(n,d)}. 
\end{align*}
Turning to the variance of $T$ and letting $\boldsymbol{\sigma} := (\sigma_1,\ldots,\sigma_d)$, the law of total variance yields
\begin{align*}
	\mV[T] ~=~ \mE\bigl[\mV[T \,|\, \boldsymbol{\sigma}]\bigr] + \mV\bigl[\mE[T \,|\, \boldsymbol{\sigma}]\bigr].
\end{align*}
The analysis of the first term is essentially the same as Proposition 3.2 of \cite{canonne2018testing}. First note that conditional on $\boldsymbol{\sigma}$, $D_1, \ldots, D_d$ are mutually independent (given $\sigma_1,\ldots,\sigma_d$, we independently generate $D_1,\ldots,D_d$), and therefore the following identity holds
\begin{align*}
	\mV[T\,|\,\boldsymbol{\sigma}] ~=~ \sum_{m \in [d]} \mathds{1}(\sigma_m \geq 4) \sigma_m^2 \mV[U(D_m)\,|\,\sigma_m].
\end{align*} 
We then use the variance bound for $U(D_m)$ in equation~(4) of \cite{canonne2018testing} and see
\begin{align*}
	\mV[T\,|\,\boldsymbol{\sigma}] ~ \lesssim ~  & \sum_{m \in [d]} \mathds{1}(\sigma_m \geq 4) \sigma_m^2 \biggl(  \frac{\|p_{X,Y\,|\,Z=m} - p_{X|Z=m}p_{Y|Z=m}\|_2^2}{\sigma_m} + \frac{1}{\sigma_m^2}  \biggr) \\[.5em]
	\lesssim ~ &  \sum_{m \in [d]} \mathds{1}(\sigma_m \geq 4) \sigma_m \|p_{X,Y\,|\,Z=m} - p_{X|Z=m}p_{Y|Z=m}\|_2^2 +   \sum_{m \in [d]} \mathds{1}(\sigma_m \geq 4) \\[.5em]
	\lesssim ~ & \mE[T\,|\,\boldsymbol{\sigma}] + \min(n, d).
\end{align*}
This reveals that $\mE\bigl[\mV[T\,|\,\boldsymbol{\sigma}]\bigr] \lesssim \mE[T] + \min(n,d)$. For the second term, the analysis is again similar to that in Proposition 3.2 of \cite{canonne2018testing}. The only difference is that $\sigma_1,\ldots,\sigma_d$ are no longer independent but negatively associated. In \cite{canonne2018testing}, the independence between $\sigma_1,\ldots,\sigma_d$ is used to show that 
\begin{align*}
	\mV\bigl[\mE[T\,|\,\boldsymbol{\sigma}]\bigr] = \mV \Biggl[ \sum_{m \in [d]} \sigma_m \mathds{1}(\sigma_m \geq 4) \delta_m^2 \Biggr] =  \sum_{m \in [d]}  \delta_m^4 \mV \bigl[ \sigma_m \mathds{1}(\sigma_m \geq 4)  \bigr].
\end{align*}
In our case, we make use of Lemma~\ref{Lemma: NA} (negative association) and prove that 
\begin{align*}
	& \mV \Biggl[ \sum_{m \in [d]} \sigma_m \mathds{1}(\sigma_m \geq 4) \delta_m^2 \Biggr]  \\[.5em]
	=~ & \sum_{m \in [d]}  \delta_m^4 \mV \bigl[ \sigma_m \mathds{1}(\sigma_m \geq 4) \bigr] + \sum_{\substack{m,m' \in [d] \\ m \neq m'}} \mathrm{Cov}\bigl\{ \sigma_m \mathds{1}(\sigma_m \geq 4), \sigma_{m'} \mathds{1}(\sigma_{m'} \geq 4) \bigr\} \\[.5em]
	\leq ~ &  \sum_{m \in [d]}  \delta_m^4 \mV \bigl[ \sigma_m \mathds{1}(\sigma_m \geq 4) \bigr].
\end{align*}
This does not change the overall conclusion as our main interest is in bounding the variance (up to constants). Additionally, Lemma~\ref{Lemma: Variance bound} yields that $\mV \bigl[ \sigma_m \mathds{1}(\sigma_m \geq 4)  \bigr] \lesssim \mE \bigl[ \sigma_m \mathds{1}(\sigma_m \geq 4)  \bigr]$ and thus 
\begin{align*}
	\mV\bigl[\mE[T \,|\, \boldsymbol{\sigma}]\bigr] \lesssim \sum_{m \in [d]}  \delta_m^4 \mE \bigl[ \sigma_m \mathds{1}(\sigma_m \geq 4)  \bigr] \lesssim \mE[T]. 
\end{align*}
Combining the previous bounds, we have $\mV[T] \lesssim \min(n,d) + \mE[T]$. Given the bounds on the expectation and the variance, the claim for $T$ follows by Chebyshev's inequality and adjusting the value of $\beta'$ and $\zeta$ as in Section 3.1.3 of \cite{canonne2018testing}.

\subsubsection{Claim for $T_W$}  \label{Section: Claim for Tw}
Next we prove the sample complexity result for $T_W$. The type II error result based on $T_W$ under Poisson sampling is given in Lemma 5.6 of \cite{canonne2018testing}. This lemma builds upon several preliminary results proved under Poisson sampling. Once we depoissonize these preliminary results, the rest of the proof remains the same. Instead of reproducing the entire proof, we only go over places where Poissonization trick is critical and reprove the corresponding results under multinomial sampling. 

\begin{itemize}
	\item The first place where Poissonization plays a role is Lemma 5.2 of \cite{canonne2018testing}. In the process of investigating the expectation of $T_W$, \cite{canonne2018testing} introduce a quantity $D = \sum_{m \in [d]} D_m$ where
	\begin{align*}
		D_m = \sigma_m \omega_m \frac{\varepsilon_m^2}{\ell_1 \ell_2} \mathds{1}(\sigma_m \geq 4) \quad \text{and} \quad \varepsilon_m = \mathrm{TV}\bigl(p_{X,Y\,|\,Z=m}, \, p_{X|Z=m}p_{Y|Z=m}\bigr).
	\end{align*}
	Then their lemma 5.2 claims that 
	\begin{align*}
		\mE[D] ~\gtrsim~ \sum_{m \in [d]} \frac{\varepsilon_m^{2}}{\ell_1 \ell_2} \min(\alpha_m \beta_m, \alpha_m^4),
	\end{align*} 
	where $\alpha_m = n p_Z(m)$ and $\beta_m =\sqrt{\min(\alpha_m, \ell_1)\min(\alpha_m, \ell_2)}$. We note that the same result holds under multinomial sampling since
	\begin{align*}
		\mE[D] ~=~ &  \sum_{m \in [d]} \frac{\varepsilon_m^2}{\ell_1 \ell_2} \mE\bigl[ \sigma_m \omega_m  \mathds{1}(\sigma_m \geq 4) \bigr] ~ \gtrsim ~ \sum_{m \in [d]} \frac{\varepsilon_m^{2}}{\ell_1 \ell_2}  \min(\alpha_m \beta_m, \alpha_m^4),
	\end{align*}	       
	where the inequality follows by Lemma~\ref{Lemma: expectation lower bound}.
	\item The second place where Poissonization is important is Lemma 5.3 of \cite{canonne2018testing}, which claims that 
	\begin{align*}
		\mV[D] ~\lesssim~ \mE[D]. 
	\end{align*}
	This result holds under multinomial sampling as well. To see this, let $\varepsilon_m' := \varepsilon_m / \sqrt{\ell_1 \ell_2}$ and observe
	\begin{align*}
		& \mV[D] ~=~  \mV\biggl[ \sum_{m \in [d]}  \sigma_m \omega_m \varepsilon_m^{\prime 2} \mathds{1}(\sigma_m \geq 4)  \biggr] \\[.5em]
		= ~ &  \sum_{m \in [d]} \mV \bigl[  \sigma_m \omega_m \varepsilon_m^{\prime 2} \mathds{1}(\sigma_m \geq 4) \bigr] + \sum_{\substack{m,m' \in [d] \\ m \neq m'}} \mathrm{Cov} \bigl\{\sigma_m \omega_m \mathds{1}(\sigma_m \geq 4),  \sigma_{m'} \omega_{m'} \mathds{1}(\sigma_{m'} \geq 4) \bigr\} \\[.5em]
		\leq ~ & \sum_{m \in [d]} \mV \bigl[  \sigma_m \omega_m \varepsilon_m^{\prime 2} \mathds{1}(\sigma_m \geq 4) \bigr],
	\end{align*}
	where the inequality uses the negative association result established in Lemma~\ref{Lemma: NA}. We then use Lemma~\ref{Lemma: Variance bound Claim 22} and see that
	\begin{align*}
		\sum_{m \in [d]} \mV \bigl[  \sigma_m \omega_m \varepsilon_m^{\prime 2} \mathds{1}(\sigma_m \geq 4) \bigr] ~ \lesssim ~ & \sum_{m \in [d]} \varepsilon_m^{\prime 4} \sqrt{\ell_1 \ell_2} \mE \bigl[  \sigma_m \omega_m \mathds{1}(\sigma_m \geq 4) \bigr] \\[.5em]
		\lesssim ~ & \sum_{m \in [d]} \varepsilon_m^{\prime 2} \mE \bigl[  \sigma_m \omega_m \mathds{1}(\sigma_m \geq 4) \bigr] = \mE[D].
	\end{align*}
	Therefore, Lemma 5.3 of \cite{canonne2018testing} holds under multinomial sampling as well. 
	\item Lemma 5.4 of \cite{canonne2018testing} provides a lower bound for the (conditional) expected value of the test statistic and an upper bound for the (conditional) variance of the test statistic, which simultaneously hold with high probability. In this result, the authors use a concentration property of $N \sim \mathrm{Poisson}(n)$ to its mean in order to connect $\min(N,d)$ with $\min(n,d)$. This step can be skipped under multinomial sampling as $n$ is fixed in this case.
\end{itemize}
Other than Lemma 5.2, Lemma 5.3 and Lemma 5.4 of \cite{canonne2018testing}, the other steps of the proof of Lemma 5.6 in \cite{canonne2018testing} are essentially the same under multinomial sampling. Hence we omit the details and finish the proof.

\subsection{Proof of Theorem~\ref{Theorem: Multinomial Sampling using Permutation}}

Our proof of Theorem~\ref{Theorem: Multinomial Sampling using Permutation} leverages the technical tools developed in \cite{kim2022minimax} and \cite{kim2021local} for controlling the randomness of permuted statistics.

\subsubsection{Claim for $T$} 
We begin with the first part that proves the sample complexity of \texttt{UCI}-test in Algorithm~\ref{Algorithm: Unweighted U-statistic} based on the test statistic $T$. By Lemma~\ref{Lemma: quantile}, the type II error of \texttt{UCI}-test is equal to 
\begin{align*}
	\mP(p_{\mathrm{perm}} > \alpha) = \mP(T \leq q_{1-\alpha}),
\end{align*}
where $q_{1-\alpha}$ is the $1-\alpha$ quantile of the empirical distribution of $T,T_1,\ldots,T_B$. Our main strategy is to show that $q_{1-\alpha}$ is upper bounded by the critical value used in Theorem~\ref{Theorem: Multinomial Sampling}. More specifically, suppose that there exists a sufficiently large $\zeta >0$ (not necessarily the same $\zeta$ in Theorem~\ref{Theorem: Multinomial Sampling}) such that 
\begin{align} \label{Eq: goal in T}
	\mP\bigl(q_{1-\alpha} \geq \zeta \sqrt{\min(n,d)}\bigr) \leq \frac{\beta}{2}. 
\end{align}
Then by denoting the event of $q_{1-\alpha} < \zeta \sqrt{\min(n,d)}$ by $\mathcal{A}$, the type II error of \texttt{UCI}-test is bounded by 
\begin{align*}
	\mP(p_{\mathrm{perm}} > \alpha) ~ \leq ~ & \mP(T \leq q_{1-\alpha}, \, \mathcal{A}) + \mP(\mathcal{A}^c) \\[.5em]
	\leq ~ & \mP\bigl(T \leq \zeta \sqrt{\min(n,d)}\bigr) + \frac{\beta}{2} \\[.5em]
	\leq ~ & \beta,
\end{align*}
where the last inequality uses the result of Theorem~\ref{Theorem: Multinomial Sampling}. Having this inequality, we only need to prove inequality~\eqref{Eq: goal in T}. To this end, we consider the permutation distribution of $T$ based on all possible local permutations. To explain it further, recall that $\Pi_{\sigma_m}$ is the set of all permutations of $[\sigma_m]$ for each $m \in [d]$, and $D_m^\pi$ is the set of the data with $Z=m$ whose $Y$ values are permuted based on $\pi \in \Pi_{\sigma_m}$. Let $\pi_{(m)}$ be a random permutation drawn uniformly from $\Pi_{\sigma_m}$ for each $m \in [d]$. Further assume that $\pi_{(1)},\ldots, \pi_{(m)}$ are mutually independent and write $\boldsymbol{\pi} = \{\pi_{(1)},\ldots, \pi_{(m)}\}$. Let $T^{\boldsymbol{\pi}}$ be the test statistic computed as $T$ but based on $\{D_1^{\pi_{(1)}},\ldots,D_d^{\pi_{(d)}}\}$. Then the permutation distribution of $T$ based on all possible local permutations is given as
\begin{align*}
	F_{\boldsymbol{\pi},T} (x) := \mE_{\pi_{(1)},\ldots,\pi_{(d)}}\bigl[  \mathds{1}\bigl(T^{\boldsymbol{\pi}} \leq x \bigr) \, \big| \, \{(X_i,Y_i,Z_i)\}_{i=1}^n \bigr]. 
\end{align*}
In other words, $F_{\boldsymbol{\pi},T}$ is the conditional distribution of $T_j$ in Algorithm~\ref{Algorithm: Unweighted U-statistic} given $\{(X_i,Y_i,Z_i)\}_{i=1}^n$. We denote the $1-\alpha$ quantile of $F_{\boldsymbol{\pi},T}$ by $q_{1-\alpha}^\dagger$. Unlike $q_{1-\alpha}$, the quantile of $F_{\boldsymbol{\pi},T}$ does not involve randomness from Monte Carlo simulation and thus is relatively easier to handle theoretically. In fact, when $B$ is sufficiently large, $q_{1-\alpha}$ and $q_{1-\alpha}^\dagger$ are close to each other by the law of large numbers. For our purpose, it is enough to show that $q_{1-\alpha}^\dagger$ is larger than $q_{1-\alpha}$ with high probability. To prove this statement, note that the Dvoretzky--Kiefer--Wolfowitz inequality guarantees that
\begin{align*}
	\sup_{x \in \mathbb{R}} \bigg| F_{\boldsymbol{\pi},T} (x) -  \frac{1}{B} \sum_{i=1}^B \mathds{1} \big( T_i \leq x \big) \bigg| \leq \sqrt{\frac{1}{2B} \log\left(\frac{4}{\beta}\right)}
\end{align*}
holds with probability at least $1-\beta/4$. Under this good event, it holds that 
\begin{align*}
	q_{1-\alpha} ~=~& \inf \bigg\{x \in \mathbb{R}: 1-\alpha \leq  \frac{1}{B+1} \sum_{i=0}^B \mathds{1} \big( T_i \leq x \big) \bigg\} \\[.5em]
	\leq ~ & \inf \bigg\{x \in \mathbb{R}: 1-\alpha \leq  \frac{1}{B+1} \sum_{i=1}^B \mathds{1} \big( T_i \leq x \big) \bigg\} \\[.5em]
	= ~ & \inf \bigg\{x \in \mathbb{R}: (1-\alpha) \frac{B+1}{B} \leq  \frac{1}{B} \sum_{i=1}^B \mathds{1} \big( T_i \leq x \big) \bigg\} \\[.5em]
	\leq ~ & \inf \bigg\{x \in \mathbb{R}: (1-\alpha) \frac{B+1}{B} + \sqrt{\frac{1}{2B} \log\left(\frac{4}{\beta}\right)} \leq  F_{\boldsymbol{\pi},T} (x) \bigg\} \\[.5em]
	\leq ~ & q_{1-\alpha/2}^\dagger, 
\end{align*}
where we utilize the assumption on $B \geq \max\{4(1-\alpha)\alpha^{-1}, \, 8\alpha^{-2}\log(4/\beta)\}$ to ensure that 
\begin{align*}
	(1-\alpha) \frac{B+1}{B} + \sqrt{\frac{1}{2B} \log\left( \frac{4}{\beta} \right)} \leq 1 - \frac{\alpha}{2}.
\end{align*}
Here constant $2$ in $\alpha/2$ is chosen for convenience and can be an arbitrary positive number greater than one. This result combined with the union bound yields that
\begin{align*}
	\mP\bigl(q_{1-\alpha} \geq \zeta \sqrt{\min(n,d)}\bigr) ~\leq~ 	\mP\bigl(q_{1-\alpha/2}^\dagger \geq \zeta \sqrt{\min(n,d)}\bigr) + \frac{\beta}{4}. 
\end{align*} 
Hence condition~\eqref{Eq: goal in T} is ensured if 
\begin{align} \label{Eq: bounding q}
	\mP\bigl(q_{1-\alpha/2}^\dagger \geq \zeta \sqrt{\min(n,d)}\bigr) \leq \frac{\beta}{4}. 
\end{align}
Indeed, the above inequality is essentially proved in the proof of Theorem 5 in \cite{kim2021local}. Strictly speaking, \cite{kim2021local} prove this under Poissonization, but this is not critical in the proof of \eqref{Eq: bounding q} and the same lines of the proof go through under multinomial sampling. This completes the proof of the first part of Theorem~\ref{Theorem: Multinomial Sampling using Permutation}. 

\subsubsection{Claim for $T_W$} 
The proof of the second part is similar to the first part. Again, we consider the permutation distribution of $T_W$ where the permutation only applies to $D_{X,Y,1},\ldots,D_{X,Y,d}$. More specifically, the permutation distribution function of $T_W$ at $x \in \mathbb{R}$ is given as the conditional expectation of $\mathds{1}(T_{j,W} \leq x)$ given everything except permutations. We denote the $1-\alpha$ quantile of the resulting permutation distribution by $q_{1-\alpha,W}^\dagger$ and the $1-\alpha$ quantile of $T_W, T_{1,W},\ldots,T_{B,W}$ by $q_{1-\alpha,W}$. By following the same logic in the first part of the proof, when $B \geq \max\{4(1-\alpha)\alpha^{-1}, \, 8\alpha^{-2}\log(4/\beta)\}$, it can be ensured that $q_{1-\alpha,W} \leq q_{1-\alpha,W}^\dagger$ with probability at least $1-\beta$. Therefore, the type II error of the permutation test in Algorithm~\ref{Algorithm: Weighted U-statistic} is bounded by
\begin{align*}
	\mP(T_W \leq q_{1-\alpha,W}) \leq \mP(T_W \leq q_{1-\alpha,W}^\dagger) + \frac{\beta}{4}.
\end{align*}  
Under Poissonization, Theorem 5 of \cite{kim2021local} shows that 
\begin{align*}
	\mP(T_W \leq q_{1-\alpha,W}^\dagger) \leq  \frac{3}{4}\beta. 
\end{align*}
In fact, the same lines of the proof of Theorem 5 in \cite{kim2021local} combined with the results in Appendix~\ref{Section: Claim for Tw} prove that the same conclusion holds under multinomial sampling. Therefore we omit the details and finish the proof.

\subsection{Proof of Proposition~\ref{Proposition: Computationally efficient formula}}
We start by making an important observation that $U^{\boldsymbol{\eta},\boldsymbol{\upsilon}}_W (D)$ is equivalent to the U-statistic of Hilbert--Schmidt Independence Criterion~\citep{gretton2005measuring} based on kernels 
\begin{align*}
	k(x,x') = \sum_{q=1}^{\ell_1} \frac{\mathds{1}(x=q)\mathds{1}(x'=q)}{\eta_q} \quad \text{and} \quad l(y,y') = \sum_{r=1}^{\ell_2} \frac{\mathds{1}(y=r)\mathds{1}(y'=r)}{\upsilon_r},
\end{align*}
for $x,x' \in [\ell_1]$ and $y,y' \in [\ell_2]$. Leveraging this observation along with the computational trick introduced in \cite{song2012feature}, $U^{\boldsymbol{\eta},\boldsymbol{\upsilon}}_W (D)$ can be computed via convenient matrix operations. Specifically, we let $\bK, \bL$ be $\sigma \times \sigma$ kernel matrices whose entries are $\bK_{i,j} = k(X_i,X_j)$ and $\bL_{i,j} = l(Y_i,Y_j)$, respectively. Define $\widetilde{\bK}$ by letting $\widetilde{\bK}_{i,j} = \bK_{i,j}$ if $i \neq j$ and $\widetilde{\bK}_{i,j}  = 0$ if $i = j$. That is, $\widetilde{\bK}$ is the same as $\bK$ except having zero diagonal components. We similarly define $\widetilde{\bL}$ by relating it to $\bL$. Denote the vector of size $\sigma$ with each entry equaling one by $\bOne$. Then equation~(5) of \cite{song2012feature} gives a computationally efficient form of $U^{\boldsymbol{\eta},\boldsymbol{\upsilon}}_W (D)$ as
\begin{align*}
	U^{\boldsymbol{\eta},\boldsymbol{\upsilon}}_W (D)~=~ \frac{1}{\sigma(\sigma-3)} \Biggl[ \tr(\widetilde{\bK} \widetilde{\bL}) + \frac{\bOne^\top \widetilde{\bK} \bOne \bOne^\top \widetilde{\bL} \bOne}{(\sigma-1)(\sigma-2)} - \frac{2}{\sigma - 2} \bOne^\top \widetilde{\bK} \widetilde{\bL} \bOne \Biggr]. 
\end{align*}
It is therefore enough to show that $A_1 = \tr(\widetilde{\bK} \widetilde{\bL})$, $A_2 = \bOne^\top \widetilde{\bK} \bOne \bOne^\top \widetilde{\bL} \bOne$ and $A_3 = \bOne^\top \widetilde{\bK} \widetilde{\bL} \bOne$. We prove these in order. 
\begin{itemize}
	\item For the claim $A_1 = \tr(\widetilde{\bK} \widetilde{\bL})$, we observe that  
	\begin{align*}
		\tr(\widetilde{\bK} \widetilde{\bL}) ~=~ &   \sum_{i=1}^n\sum_{j=1}^n \bK_{i,j} \bL_{i,j} \mathds{1}(i \neq j) \\[.5em]
		= ~ & \sum_{i=1}^n \sum_{j=1}^n \Biggl[ \sum_{q=1}^{\ell_1} \frac{\mathds{1}(X_i=q)\mathds{1}(X_j=q)}{\eta_q} \Biggr] \Biggl[ \sum_{r=1}^{\ell_2} \frac{\mathds{1}(Y_i=r)\mathds{1}(Y_j=r)}{\upsilon_r} \Biggr]  \mathds{1}(i \neq j) \\[.5em] 
		= ~ & \sum_{q=1}^{\ell_1} \sum_{r=1}^{\ell_2} \sum_{i=1}^n \sum_{j=1}^n \frac{\mathds{1}(X_i=q) \mathds{1}(X_j=q) \mathds{1}(Y_i = r) \mathds{1}(Y_j=r)}{\eta_q \upsilon_r} \mathds{1}(i \neq j) \\[.5em]
		= ~ &  \sum_{q=1}^{\ell_1}  \sum_{r=1}^{\ell_2} \Biggl[ \frac{o_{qr}^2 - o_{qr}}{\eta_q \upsilon_r} \Biggr] = A_1.
	\end{align*}
	\item To prove the second claim $A_2 = \bOne^\top \widetilde{\bK} \bOne \bOne^\top \widetilde{\bL} \bOne$, notice that 
	\begin{align*}
		\bOne^\top \widetilde{\bK} \bOne ~=~ & \sum_{i=1}^n\sum_{j=1}^n \sum_{q=1}^{\ell_1} \frac{\mathds{1}(X_i= q) \mathds{1}(X_j=q)}{\eta_q} \mathds{1}(i \neq j) \\[.5em]
		= ~ & \sum_{q=1}^{\ell_1} \Biggl[ \frac{o_{q +}^2 - o_{q +}}{\eta_q} \Biggr].
	\end{align*}
	Similarly, we have 
	\begin{align*}
		\bOne^\top \widetilde{\bL} \bOne ~=~ & \sum_{r=1}^{\ell_2} \Biggl[ \frac{o_{+ r}^2 - o_{+ r}}{\upsilon_r} \Biggr].
	\end{align*}
	Combining these two proves the identity $A_2 = \bOne^\top \widetilde{\bK} \bOne \bOne^\top \widetilde{\bL} \bOne$.
	\item Finally, the claim $A_3 = \bOne^\top \widetilde{\bK} \widetilde{\bL} \bOne$ can be proved by following a series of identities:
	\begin{align*}
		\bOne^\top \widetilde{\bK} \widetilde{\bL} \bOne ~=~ & \sum_{i=1}^n \Biggl[ \sum_{j=1}^n \bK_{i,j}  \mathds{1}(i \neq j) \Biggr] \Biggl[ \sum_{j'=1}^n \bL_{i,j'}  \mathds{1}(i \neq j')\Biggr]  \\[.5em]
		= ~ &\sum_{i=1}^n \Biggl[ \sum_{j=1}^n \sum_{q=1}^{\ell_1} \frac{\mathds{1}(X_i=q)\mathds{1}(X_j=q)}{\eta_q} \mathds{1}(i \neq j) \Biggr]  \Biggl[ \sum_{j'=1}^n \sum_{r=1}^{\ell_2} \frac{\mathds{1}(Y_i=r)\mathds{1}(Y_{j'}=r)}{\upsilon_r} \mathds{1}(i \neq j') \Biggr] \\[.5em]
		= ~ & \sum_{q=1}^{\ell_1} \sum_{r=1}^{\ell_2} \sum_{i=1}^n \Biggl[ \sum_{j=1}^n \frac{\mathds{1}(X_i=q) \mathds{1}(X_j=q)}{\eta_q} \mathds{1}(i \neq j) \Biggr]  \Biggl[ \sum_{j'=1}^n \frac{\mathds{1}(Y_i=r) \mathds{1}(Y_{j'}=r)}{\upsilon_r} \mathds{1}(i \neq j') \Biggr] \\[.5em]
		= ~ &  \sum_{q=1}^{\ell_1} \sum_{r=1}^{\ell_2} \frac{o_{qr} o_{q +} o_{+ r}}{\eta_q \upsilon_r} -  \sum_{q=1}^{\ell_1} \sum_{r=1}^{\ell_2} \frac{o_{qr} o_{q+}}{\eta_q \upsilon_r} -  \sum_{q=1}^{\ell_1} \sum_{r=1}^{\ell_2} \frac{o_{qr} o_{+ r}}{\eta_q \upsilon_r} +  \sum_{q=1}^{\ell_1} \sum_{r=1}^{\ell_2} \frac{o_{qr}}{\eta_q \upsilon_r} \\[.5em]
		= ~ &  \sum_{q=1}^{\ell_1} \sum_{r=1}^{\ell_2} \frac{o_{qr}(o_{q+}o_{+r} - o_{q+} - o_{+r} + 1)}{\eta_q \upsilon_r} \\[.5em]
		= ~ & A_3.
	\end{align*}
\end{itemize}
This completes the proof of Proposition~\ref{Proposition: Computationally efficient formula}.

\subsection{Proof of Proposition~\ref{Proposition: sub-optimality}} \label{Proof: Proposition: sub-optimality}
Our strategy is to find a sequence of distributions under the given conditions for which the asymptotic power of $\chi^2$- and $G$-tests are zero. This result implies the desired claim since we are concerned with the worse case power. Let us start by forming a sequence of discrete distributions on $[2] \times [2] \times [d]$. 

\vskip 1em 

\noindent \textbf{Worst case scenario.} Consider the setting where $p_{X,Y|Z}(1,1\,|\,z)=p_{X,Y|Z}(2,2\,|\,z) = 1/2$ and $p_{X,Y|Z}(1,2\,|\,z)=p_{X,Y|Z}(2,1\,|\,z) = 0$ for all $z \in [d]$. For each $n$, set $p_Z(1) = (1 - 1/n)^{1/n}$ so that the probability of having all of $Z_1,\ldots,Z_n$ equal to one is $1 - 1/n$. We then set $p_Z(2) = \ldots = p_Z(d) = \bigl(1 - p_Z(1)\bigr)/(d-1)$. Notice that by construction, $p_Z$ takes non-zero values over $[d]$ and thus the dimension of $Z$ is $d$. We then let $p = p_{X,Y|Z} p_Z$ and claim that 
\begin{align} \label{Eq: claim of lower bound}
	\inf_{q \in \mathcal{P}_0} \|p - q\|_1 \geq \frac{1}{4}.
\end{align}
Note that for any $q = q_{X|Z}q_{Y|Z}q_Z \in \mathcal{P}_0$,
\begin{align*}
	\|p - q\|_1 ~=~ & \|p_{X,Y|Z}p_Z - q_{X|Z}q_{Y|Z}q_Z\|_1 \\[.5em]
	\overset{\mathrm{(i)}}{\geq} ~ & \|p_{X,Y|Z}p_Z - q_{X|Z}q_{Y|Z}p_Z\|_1 -  \|q_{X|Z}q_{Y|Z}p_Z - q_{X|Z}q_{Y|Z}q_Z\|_1 \\[.5em]
	= ~ & \sum_{z \in [d]} p_Z(z) \cdot \|p_{X,Y|Z=z} - q_{X|Z=z}q_{Y|Z=z}\|_1 -  \|p_Z - q_Z\|_1 \\[.5em]
	\overset{\mathrm{(ii)}}{\geq} ~ & \sum_{z \in [d]} p_Z(z) \cdot \|p_{X,Y|Z=z} - q_{X|Z=z}q_{Y|Z=z}\|_1   - 	\|p - q\|_1,
\end{align*}
where step~(i) uses the triangle inequality and step~(ii) uses the following inequality (which can be proved using the triangle inequality):
\begin{align*}
	\|p_Z - q_Z\|_1 = \sum_{z \in [d]} | p_Z(z) - q_Z(z) | = \sum_{z \in [d]} \Bigg| \sum_{x \in [\ell_1], y \in [\ell_2] }p(x,y,z) - q(x,y,z) \Biggr| \leq \|p - q\|_1. 
\end{align*}
Therefore we can conclude that
\begin{align*}
	\inf_{q \in \mathcal{P}_0} \|p - q\|_1 \geq  \frac{1}{2} \inf_{q \in \mathcal{P}_0} \Biggl\{ \sum_{z \in [d]} p_Z(z) \cdot \|p_{X,Y|Z=z} - q_{X|Z=z}q_{Y|Z=z}\|_1 \Biggr\}.
\end{align*}
In order to further lower bound the above quantity, we first prove a preliminary result. Let $\mathcal{Q}_{0}$ be the collection of joint distributions of $(X,Y)$ on $[2] \times [2]$ such that for any $q_{X,Y} \in \mathcal{Q}_{0}$, $q_{X,Y} = q_X q_Y$, i.e.,~$X \, \indep \, Y$. Then for a joint distribution $p_{X,Y}$ such that  $p_{X,Y}(1,1)=p_{X,Y}(2,2) = 1/2$ and $p_{X,Y}(1,2)=p_{X,Y}(2,1) = 0$, we have 
\begin{align*}
	\inf_{q \in \mathcal{Q}_0} \| p_{X,Y} - q \|_1 \geq \inf_{q \in \mathcal{Q}_0} \| p_{X,Y} - q \|_2 = \frac{1}{2},
\end{align*}
where the last equality can be shown by solving the following optimization problem:
\begin{align*}
	\inf_{a,b \in [0,1]} \sqrt{\Bigl( \frac{1}{2} - ab \Bigr)^2 + a^2(1-b)^2 + (1-a)^2 b^2 + \Bigl( \frac{1}{2} - (1-a)(1-b) \Bigr)^2} = \frac{1}{2},
\end{align*}
where the infimum is achieved at $a = b = 1/2$. Consequently, given that $(X,Y) \, \indep \, Z$ in our construction, 
\begin{align*}
	\inf_{q \in \mathcal{P}_0} \|p - q\|_1 ~\geq~& \frac{1}{2} \inf_{q \in \mathcal{P}_0} \Biggl\{ \sum_{z \in [d]} p_Z(z) \cdot \|p_{X,Y|Z=z} - q_{X|Z=z}q_{Y|Z=z}\|_1 \Biggr\} \\[.5em]
	\geq ~ & \frac{1}{4},
\end{align*} 
which proves claim~\eqref{Eq: claim of lower bound}. 

\vskip 1em

\noindent \textbf{Main proof.} Now let $\mathcal{A}_n$ be an event that all of $Z_1,\ldots,Z_n$ are equal to one, and $\mathcal{A}^c$ be its complement. Starting with $\chi^2$-test, its worst case power is bounded by 
\begin{equation} \label{Eq: upper bound}
	\begin{aligned}
		\inf_{p \in \mathcal{P}_1(\varepsilon)}\mP_{p}(\chi^2 > q_{1-\alpha,d}) ~\leq~ & \mP_{p^\ast}(\chi^2 > q_{1-\alpha,d}) \\[.5em]
		\leq ~& \mP_{p^\ast}(\chi^2 > q_{1-\alpha,d} \, | \,  \mathcal{A}_n) \mP_{p^\ast}(\mathcal{A}_n) +  \mP_{p^\ast}(\mathcal{A}_n^c) \\[.5em] 
		\leq ~& \mP_{p^\ast}(\chi^2 > q_{1-\alpha,d} \,  | \,  \mathcal{A}_n) + \frac{1}{n}.
	\end{aligned}
\end{equation}
Conditional on $\mathcal{A}_n$, $\chi^2/n$ has the same distribution as 
\begin{align*}
	W^\star_n := \frac{(W-W^2/n)^2}{W^2} + \frac{2\bigl(W^2(n-W)^2/n^2\bigr)}{W(n-W)} + \frac{\bigl(n-W-(n-W)^2/n\bigr)^2}{(n-W)^2},
\end{align*}
where $W \sim \mathrm{Binomial}(n,1/2)$. Using this notation, we have the identity
\begin{align*}
	\mP_{p^\ast}(\chi^2 > q_{1-\alpha,d} \,  | \,  \mathcal{A}_n) ~=~ \mP\biggl( \frac{n}{d} \cdot W^\star_n- 1 > \frac{1}{\sqrt{d}} \cdot \frac{q_{1-\alpha,d} - d}{\sqrt{d}} \biggr).
\end{align*}
Under our choice of $d = n \times r_n$ where $r_n \rightarrow \infty$, we apply the law of large number and continuous mapping theorem and see that
\begin{align*}
	 \frac{n}{d} \cdot W^\star_n- 1 \convP -1,
\end{align*}
where $X_n \convP X$ means convergence of $X_n$ to $X$ in probability. In addition, by the central limit theorem along the fact that convergence in distribution implies convergence of the quantile function at all continuity points, we have
\begin{align*}
	\frac{q_{1-\alpha} - d}{\sqrt{2d}} \rightarrow z_{1-\alpha} \quad \text{as $d \rightarrow \infty$,}
\end{align*}
where $z_{1-\alpha}$ is the $1-\alpha$ quantile of $N(0,1)$. Thus
\begin{align*}
	 & \frac{1}{\sqrt{d}} \cdot \frac{q_{1-\alpha,d} - d}{\sqrt{d}} \rightarrow 0 \quad \text{as $n \rightarrow \infty$.}
\end{align*}
Therefore by Slutsky's theorem, we have
\begin{align*}
	\mP_{p^\ast}(\chi^2 > q_{1-\alpha,d} \,  | \,  \mathcal{A}_n) = \mP\biggl( \frac{n}{d} \cdot W^\star_n- 1 > \frac{1}{\sqrt{d}} \cdot \frac{q_{1-\alpha,d} - d}{\sqrt{d}} \biggr) \rightarrow 0. 
\end{align*}
This in conjunction with the upper bound in \eqref{Eq: upper bound} implies that 
\begin{align*}
	\lim_{n \rightarrow \infty} \inf_{p \in \mathcal{P}_1(\varepsilon)}\mP_{p}(\chi^2 > q_{1-\alpha,d}) = 0. 
\end{align*}
The proof of $G$-test is essentially the same. The only difference is that conditional on $\mathcal{A}_n$, $G/n$ has the same distribution as 
\begin{align*}
	V_n^\star := 2 \biggl[ \frac{W}{n} \log \biggl( \frac{n}{W} \biggr) + \biggl( 1 - \frac{W}{n} \biggr) \! \log \biggl( \frac{1}{1 - \frac{W}{n}} \biggr) \biggr],
\end{align*}
from which it holds that 
\begin{align*}
	\mP_{p^\ast}(G > q_{1-\alpha,d} \,  | \,  \mathcal{A}_n) ~=~ \mP\biggl( \frac{n}{d} \cdot V^\star_n - 1 > \frac{1}{\sqrt{d}} \cdot \frac{q_{1-\alpha,d} - d}{\sqrt{d}} \biggr).
\end{align*}
Again, it can be seen that $\frac{n}{d} \cdot V^\star_n - 1 \convP -1$ and $\frac{1}{\sqrt{d}} \cdot \frac{q_{1-\alpha,d} - d}{\sqrt{d}} \rightarrow 0$. Following the same steps as in the case of $\chi^2$-test, we conclude 
\begin{align*}
	\lim_{n \rightarrow \infty} \inf_{p \in \mathcal{P}_1(\varepsilon)}\mP_{p}(G > q_{1-\alpha,d}) = 0. 
\end{align*}
This completes the proof of Proposition~\ref{Proposition: sub-optimality}.

\subsection{Proof of Lemma~\ref{Lemma: lower bound}}
We start by computing the expectation of $X \mathds{1}(X \geq 4)$ explicitly as
\begin{align*}
	\mE[X \mathds{1}(X \geq 4)] ~=~ & \mE[X] - \mE[X \mathds{1}(X \leq 3)] \\[.5em]
	= ~ & np - np(1-p)^{n-1} - n(n-1)p^2(1-p)^{n-2}-\frac{n(n-1)(n-2)}{2}p^3(1-p)^{n-3}.
\end{align*}
Given this formula, we consider two cases: (i) $np > 1$ and (ii) $np \leq 1$, separately. Throughout the proof, we assume that $p>0$ since the claim is trivial when $p=0$.

\vskip 1em 

\noindent \textbf{Case (i) $np > 1$.} First assume that $np > 1$, equivalently $p > 1/n$. At a high-level, when $np$ is large, there is only a small possibility that $\mathds{1}(X < 4 ) = 1$. In this case, $X \mathds{1}(X \geq 4)$ behaves like $X$ and we expect that $\mE[X \mathds{1}(X \geq 4)] \approx \mE[X] = np$. We now make this statement more precise. To this end, let us define 
\begin{align*}
	r_n(p) ~ := ~ & \frac{\mE[X \mathds{1}(X \geq 4)]}{np}  \\[.5em]
	= ~ &  1 - (1-p)^{n-1} - (n-1)p(1-p)^{n-2} - \frac{(n-1)(n-2)}{2}p^2(1-p)^{n-3}.
\end{align*}
We will show that $r_n(p)$ is an increasing function of $p$ for $p > 1/n$ and thus the minimum is achieved at $r_n(1/n)$. To show this, we calculate its derivative with respect to $p$, which can be seen to be non-negative for all $n \geq 4$ as 
\begin{align*}
	r_n'(p) = \frac{(n-1)(n-2)(n-3)(1-p)^np^2}{2(1-p)^4} \geq 0 \quad \text{for all $n\geq 4$.}
\end{align*}
Therefore, the minimum of $r_n(p)$ is achieved at $p = 1/n$ as a function of $p$, which yields 
\begin{align*}
	r_n(1/n) ~=~ & 1 - \left(1 - \frac{1}{n}\right)^{n-1} - \frac{n-1}{n}\left(1 - \frac{1}{n}\right)^{n-2} - \frac{(n-1)(n-2)}{2n^2} \left(1 - \frac{1}{n}\right)^{n-3} \\[.5em]
	= ~ & \biggl( 1 - \frac{1}{n} \biggr)^{n} \biggl\{\frac{n(5n-6)}{2(n-1)^2}\biggr\} \\[.5em]
	\overset{(\ast)}{\geq} ~ & e^{-1} \biggl(\frac{n-1}{n}\biggr) \biggl\{\frac{n(5n-6)}{2(n-1)^2} \biggr\} = e^{-1} \biggl\{ \frac{5}{2} - \frac{1}{2(n-1)}\biggr\},
\end{align*}		
where step~($\ast$) follows by the inequality
\begin{align*}
	\left( 1 + \frac{x}{n} \right)^n \geq e^x \left( 1 - \frac{x^2}{n}\right) \quad\text{for all $n \geq 1$ and $|x| \leq n$}.
\end{align*}
From this, we conclude that $r_n(1/n) \geq e^{-1} \times \frac{7}{3} \approx 0.858$ for all $n \geq 4$ and thus
\begin{align*}
	\mE[X\mathds{1}(X \geq 4)] \geq \gamma np \quad \text{if $np > 1$.}
\end{align*}

\vskip 1em 

\noindent \textbf{Case (ii) $np \leq 1$.} Next assume $np \leq 1$, i.e.,~$p \leq 1/n$. Then $\min\{np, (np)^4\} = (np)^4$ and the ratio function $r_n(p)$ is computed as
\begin{align*}
	r_n(p) ~ :=~ & \frac{\mE[X \mathds{1}(X \geq 4)]}{n^4p^4} \\[.5em] 
	= ~ & \frac{1}{(np)^3} - \frac{(1-p)^{n-1}}{(np)^3} - \frac{(n-1)}{n^3p^2}(1-p)^{n-2} - \frac{(n-1)(n-2)}{2n^3p} (1-p)^{n-3}.
\end{align*}
Given this formula, we aim to prove that $r_n(p)$ is a decreasing function of $p$ for $p \leq 1/n$ and hence the minimum is achieved at $r_n(1/n) \geq \gamma$ as shown earlier. To prove this, by assuming $p \neq 0$, it can be seen that the derivative of $r_n(p)$ as a function of $p$ is given by
\begin{align*}
	r_n'(p) ~=~\frac{(1-p)^n\{(n-4)p[(n-3)p((n-2)p+3)+6]+6 \} - 6(1-p)^4}{2n^3(1-p)^4p^4}. 
\end{align*}
Our goal is to show that $r_n'(p) \leq 0$ for any $p \in (0,1/n]$, which is equivalent to the conditions
\begin{align} \nonumber
	& r_n'(p) \leq 0 \\[.5em]   \nonumber
	\overset{\mathrm{iff}}{\Longleftrightarrow} ~ & 6 (1-p)^n - 6(1-p)^4+ (1-p)^n(n-4) p\{(n^2-5n+6)p^2 +3 (n-3) p + 6\} \leq 0 \\[.5em] \label{Eq: intermediate condition}
	\overset{\mathrm{iff}}{\Longleftrightarrow} ~ & 6 - 6(1-p)^{4-n}  + (n-4) \{(n-2)(n-3)p^3 +3 (n-3) p^2 + 6p \} \leq 0.
\end{align}
To this end, we make a key (and nontrivial) observation that the function $h_n(p) := 6(1-p)^{4-n}$ has a Taylor expansion around $0$ as
\begin{align*}
	 h_n(p) ~=~ & 6 + 6(n-4)p + 3(n-4)(n-3)p^2 + (n-4)(n-3)(n-2)p^3 \\[.5em]
	 & + \underbrace{\frac{(n-4)(n-3)(n-2)(n-1)}{(1-\xi)^n} p^4}_{\mathrm{remainder}},
\end{align*}
for some $\xi \in (0, 1/n)$. Importantly, the remainder term is strictly positive, which concludes that the inequality~(\ref{Eq: intermediate condition}) holds as well as $r_n'(p) \leq 0$ for all $p \in (0,1/n]$. This means that $r_n(p)$ is a decreasing function and $r_n(1/n) \leq r_n(p)$ for $p \in (0,1/n]$. From the previous calculation in case~(i), we know that $r_n(1/n) \geq e^{-1} \times \frac{7}{3}$ for all $n \geq 4$ (when $p=1/n$, we have $np = (np)^4$) and thereby prove the claim.

\subsection{Proof of Lemma~\ref{Lemma: Variance bound}} 
We divide the proof into two parts depending on whether (i)~$np > 1$ or (ii)~$np \leq 1$. Throughout, we assume that $p>0$ since the claim is trivial when $p=0$. 

\vskip 1em

\noindent \textbf{Case (i) $np > 1$.} The statement is relatively easier to prove when $np > 1$. Again, the intuition is that when $np$ is large, $X \mathds{1}(X \geq 4)$ behaves like $X$ and thus $\mathrm{Var}[X \mathds{1}(X \geq 4)] \approx \mathrm{Var}[X] = np(1-p)$ and $\mE[X \mathds{1}(X \geq 4)] \approx \mE[X] = np$. Thus the expectation is larger than the variance. To make this statement more precise, note that the variance can be bounded by
\begin{align*}
	\mathrm{Var}[X \mathds{1}(X \geq 4)]  ~=~ & \mE[X^2] - \mE[X^2 \mathds{1}(X \leq 3)] - \{ \mE[X] - \mE[X  \mathds{1}(X \geq 3)] \}^2 \\[.5em]
	\leq ~ & \mE[X^2] - \{\mE[X]\}^2 + 2\mE[X] \mE[X \mathds{1}(X \geq 3)] \\[.5em]
	\leq ~ & np(1-p) + 6 np \mP(X \geq 3) \\[.5em]
	\leq ~ & 7np.
\end{align*}
From Lemma~\ref{Lemma: lower bound}, we know that $\mE[X \mathds{1}(X \geq 4)] \ge C np$ for some $C>0$ when $np > 1$. Therefore, by taking $\gamma$ larger than, for example, $7/0.85$, the desired result follows. 

\vskip 1em

\noindent \textbf{Case (ii) $np \leq 1$. } Now we assume that $np \leq 1$, i.e.,~$p \leq 1/n$. In this case, we cannot ignore the effect of the indicator function $\mathds{1}(X \geq 4)$ as in the case of $np >1$ and thus it requires a more delicate analysis. To simplify the problem a bit, we first upper bound the ratio between the variance and the expectation as
\begin{align*}
	\frac{\mathrm{Var}[X \mathds{1}(X \geq 4)]}{\mE[X \mathds{1}(X \geq 4)]} \lesssim \frac{\mE[X^2 \mathds{1}(X \geq 4)]}{(np)^4},
\end{align*}
where we use the fact that $\mathrm{Var}[X \mathds{1}(X \geq 4)] \leq \mE[X^2 \mathds{1}(X \geq 4)]$ and $\mE[X \mathds{1}(X \geq 4)] \gtrsim (np)^4$ when $np \leq 1$ from Lemma~\ref{Lemma: lower bound}. Our goal is to show that 
\begin{align*}
	r_n(p):=\frac{\mE[X^2 \mathds{1}(X \geq 4)]}{(np)^4}
\end{align*}
is upper bounded by some positive constant for all $n \geq 4$ and $p \leq 1/n$. In particular, we will show that $r_n(p)$ is a decreasing function in $p$ for all $n \geq 4$ and
\begin{align*}
	\lim_{p \rightarrow 0} r_n(p) = \frac{2(n^3-6n^2+11n - 6)}{3n^3} < 0.67 \quad \text{for all $n \geq 4$.}
\end{align*}
To verify these claims, note that the second moment of $X \mathds{1}(X \geq 4)$:
\begin{align*}
	& \mE[X^2 \mathds{1}(X \geq 4)] ~=~  \mE[X^2] - \mE[X^2 \mathds{1}(X \leq 3)] \\[.5em]
	= ~ & np(1-p) + (np)^2  - np(1-p)^{n-1} - 2n(n-1)p^2(1-p)^{n-2} - \frac{3}{2} n(n-1)(n-2)p^3(1-p)^{n-3},
\end{align*}
and therefore $r_n(p)$ can be written as
\begin{align*}
	r_n(p) ~=~ & \frac{1-p}{(np)^3} + \frac{1}{(np)^2} - \frac{(1-p)^{n-1}}{(np)^3} - \frac{2(n-1)(1-p)^{n-2}}{n^3p^2} - \frac{3(n-1)(n-2)(1-p)^{n-3}}{2n^3p}.
\end{align*}
To show that $r_n(p)$ is a decreasing function in $p$, consider its derivative with respect to $p$:
\begin{align*}
	r_n'(p) ~=~ & \dfrac{3\left(n-3\right)\left(n-2\right)\left(n-1\right)\left(1-p\right)^{n-4}}{2n^3p}+\dfrac{7\left(n-2\right)\left(n-1\right)\left(1-p\right)^{n-3}}{2n^3p^2}+\dfrac{5\left(n-1\right)\left(1-p\right)^{n-2}}{n^3p^3} \\[.5em]
	& +\dfrac{3\left(1-p\right)^{n-1}}{n^3p^4} -\dfrac{2}{n^2p^3} - \frac{1}{n^3p^3}-\dfrac{3(1-p)}{n^3p^4}.
\end{align*}
We would like to show that $r_n'(p)$ is less than or equal to $0$ for all $0 < p \leq 1/n$ and $n \geq 4$, which is equivalent to verifying the following slightly simplified conditions:
\begin{align*}
	& r_n'(p) \leq 0 \\[.5em]
	\overset{\mathrm{iff}}{\Longleftrightarrow} \quad & \frac{3}{2}\left(n-3\right)\left(n-2\right)\left(n-1\right)np^3\left(1-p\right)^{n-4}+\frac{7}{2}\left(n-2\right)\left(n-1\right)np^2\left(1-p\right)^{n-3} \\[.5em]
	&+5\left(n-1\right)np\left(1-p\right)^{n-2} + 3n \left(1-p\right)^{n-1} - 2n^2p - np- 3n(1-p) ~ \leq ~ 0\\[.5em] 
	\overset{\mathrm{iff}}{\Longleftrightarrow} \quad & 3\left(n-3\right)\left(n-2\right)\left(n-1\right)np^3\left(1-p\right)^{n-4} + 7\left(n-2\right)\left(n-1\right)np^2\left(1-p\right)^{n-3}  \\[.5em]
	& + 10\left(n-1\right)np\left(1-p\right)^{n-2}  + 6n \left(1-p\right)^{n-1} - 4n^2p - 2np-6n(1-p) ~ \leq ~ 0 \\[.5em]
	\overset{\mathrm{iff}}{\Longleftrightarrow} \quad & 3\left(n-3\right)\left(n-2\right)\left(n-1\right)np^4\left(1-p\right)^{n-4} + 7\left(n-2\right)\left(n-1\right)np^3\left(1-p\right)^{n-3} \\[.5em]
	& + 10\left(n-1\right)np^2\left(1-p\right)^{n-2} + 6n p\left(1-p\right)^{n-1} - 4n^2p^2 - 2np^2 - 6np(1-p) ~ \leq ~ 0.
\end{align*}
Focusing on the last equivalent condition, we first define 
\begin{align*}
	f_n(p) ~:=~ &  3\left(n-3\right)\left(n-2\right)\left(n-1\right)np^4\left(1-p\right)^{n-4} + 7\left(n-2\right)\left(n-1\right)np^3\left(1-p\right)^{n-3} \\[.5em]
	& + 10\left(n-1\right)np^2\left(1-p\right)^{n-2} + 6n p\left(1-p\right)^{n-1}.
\end{align*}
Then the last condition simply becomes 
\begin{align} \label{Eq: simple claim}
	f_n(p) \leq 6np + 4n(n-1)p^2. 
\end{align}
Somewhat surprisingly, it can be seen that $6np + 4n(n-1)p^2$ is the second-order Taylor expansion of $f_n(p)$ near $p=0$. In particular, we have 
\begin{align*}
	f_n(p) ~=~ &  f_n(0) + f_n'(0) p + \frac{f_n''(0)}{2} p^2 + R_n \\[.5em]
	= ~ & 6np + 4n(n-1)p^2 + R_n,
\end{align*}
where 
\begin{align*}
	R_n = \frac{f_n^{(3)}(\zeta)}{6} p^3 \quad \text{for some $\zeta \in (0,1/n)$.}
\end{align*}
Therefore claim~(\ref{Eq: simple claim}) follows once we show that $R_n \leq 0$ for all $n \geq 4$ and $p \leq 1/n$. Note that the third derivative of $f_n(p)$ is computed as
\begin{align*}
	& f_n^{(3)}(p) \\[.5em]
	= ~ & \dfrac{n(n-1)(n-2)(n-3)(n-4)(1-p)^n p^2  \{ n(3n-4)p^2 + (35-29n)p+55\}}{(p-1)^7}
\end{align*}
and since the denominator $(p-1)^7$ is negative, it holds that $f_n^{(3)}(p) \leq 0$ once 
\begin{align*}
	k_n(p) := n(3n-4)p^2 + (35-29n)p+55
\end{align*}
is positive for $p \in (0,1/n]$. Note that $k_n(p)$ is a quadratic function of $p$ and its global minimum is achieved at $p^\ast_n = \frac{29n - 35}{2n(3n-4)}$. Moreover,
\begin{align*}
	p^\ast_n - \frac{1}{n} = \frac{23n - 27}{2n(3n-4)} > 0 \quad \text{for $n \geq 4$,}
\end{align*}
which means that $k_n(p)$ is decreasing for all $p \leq \frac{1}{n}$. On the other hand, $k_n(1/n)$ is computed as
\begin{align*}
	k_n(1/n) = 29 + \frac{23}{n} > 0 \quad \text{for all $n \geq 4$.}
\end{align*} 
Therefore $k_n(p)$ is positive on $p \in (0,1/n]$ for all $n \geq 4$, which in turn implies $R_n \leq 0$. In conclusion, the claim in \eqref{Eq: simple claim} holds for all $n \geq 4$ and therefore we complete the proof of Lemma~\ref{Lemma: Variance bound}.

\subsection{Proof of Lemma~\ref{Lemma: Variance bound Claim 22}} 
In contrast to a Poisson random variable that is characterized by a single parameter $\lambda$, a binomial random variable depends on two parameters $n$ and $p$. As in the case of Lemma~\ref{Lemma: Variance bound}, this makes the proof significantly more challenging than the proof for a Poisson random variable. Indeed, we need to take a totally different approach than \cite{canonne2018testing} and \cite{kim2022comments} to prove the given statement. The key ingredient of our proof is the Efron--Stein inequality recalled in Lemma~\ref{Lemma: Efron-Stein}. 

Let $X_1,\ldots,X_n \overset{\mathrm{i.i.d.}}{\sim} \mathrm{Binomial}(1,p)$. Then by letting $X = \sum_{i=1}^n X_i$, we can represent the quantity of interest as a function of i.i.d.~random variables as
\begin{align*}
	f(X_1,\ldots,X_n) = X\sqrt{\min\{X,a\}\min\{X,b\}} \mathds{1}(X\geq 4).
\end{align*}
For simplicity, we denote $f(X_1,\ldots,X_n)$ by $f$. Then the Efron--Stein inequality (Lemma~\ref{Lemma: Efron-Stein}) yields that 
\begin{align*}
	\mathrm{Var}\bigl[ X\sqrt{\min\{X,a\}\min\{X,b\}} \mathds{1}(X\geq 4) \bigr]  \leq \frac{1}{2} \sum_{i=1}^n \mE\big[\big(f - f^{i}\big)^2\big],
\end{align*}
where we define $f^{i} = f(X_1,\ldots,X_{i-1},X_i',X_{i+1},\ldots,X_n)$ and $X_i' \sim \mathrm{Binomial}(1,p)$ independent of $X_1,\ldots,X_n$. By symmetry, it holds that 
\begin{align*}
	\mathrm{Var}\bigl[ X\sqrt{\min\{X,a\}\min\{X,b\}} \mathds{1}(X\geq 4) \bigr]  \leq \frac{n}{2} \mE\big[\big(f - f^{1}\big)^2\big],
\end{align*}
and we let $Y = \sum_{i=2}^n X_i$ for simplicity so that $X = X_1 + Y$. By the law of total expectation, 
\begin{align*}
	n \mE\big[\big(f - f^{1}\big)^2\big] ~=~ & n \mE_Y \big[ \mE_{X_1,X_1'} \big[(f - f^1)^2 \,|\, Y \big] \big] \\[.5em]
	= ~ & \underbrace{n \mE_Y \big[ \mE_{X_1,X_1'} \big[(f - f^1)^2 \,|\, Y \big] \mathds{1}(Y < 4) \big]}_{(\mathrm{I})} +  \underbrace{n \mE_Y \big[ \mE_{X_1,X_1'} \big[(f - f^1)^2 \,|\, Y \big]\mathds{1}(Y \geq 4) \big]}_{(\mathrm{II})}. 
\end{align*}
We analyze two terms (I) and (II) separately.

\vskip 1em

\noindent \textbf{Analysis of (I).} We first claim that 
\begin{align*}
	\mE_Y \big[ \mE_{X_1,X_1'} \big[(f - f^1)^2 \,|\, Y \big] \mathds{1}(Y < 4) \big] = \mE_Y \big[ \mE_{X_1,X_1'} \big[(f - f^1)^2 \,|\, Y \big] \mathds{1}(Y = 3) \big].
\end{align*}
To see this, we unpack the expression of $f$ and $f^1$ as
\begin{align*}
	& (f - f^1)^2 \mathds{1}(Y < 4) \\[.5em]
	=~ &  \Big\{(X_1+ Y)\sqrt{\min\{(X_1+Y),a\}\min\{(X_1+Y),b\}} \mathds{1}(X_1+Y\geq 4) \\[.5em]
	& - (X_1'+Y)\sqrt{\min\{(X_1'+Y),a\}\min\{(X_1'+Y),b\}} \mathds{1}(X_1'+ Y\geq 4)  \Big\}^2  \mathds{1}(Y < 4). 
\end{align*}
When $Y \in \{0,1,2\}$, the indicator functions $\mathds{1}(X_1+Y\geq 4)$ and $\mathds{1}(X_1+Y\geq 4)$ are equal to zero since $X_1$ and $X_1'$ are either zero or one. Therefore, we only need to focus on the case of $Y = 3$. Conditional on $Y=3$, we make use of the fact that $X_1$ and $X_1'$ are independent and show that the expectation of $(f - f^1)^2 \mathds{1}(Y = 3)$ is 
\begin{align*}
	2p(1-p) (Y+1) \sqrt{\min\{(Y+1),a\}\min\{(Y+1),b\}} \mathds{1}(Y=3). 
\end{align*}
Now by taking the expectation over $Y \sim \mathrm{Binomial}(n-1,p)$,
\begin{align*}
	\mE_Y \Big[ \mE_{X_1,X_1'} \big[(f - f^1)^2 \,|\, Y \big] \mathds{1}(Y < 4) \Big] ~ = ~ & 8p(1-p)   \sqrt{\min\{4,a\}\min\{4,b\}} \binom{n-1}{3}p^3(1-p)^{n-4}.
\end{align*} 
Therefore by using the condition that $n \geq 4$ and the inequality $(1+x) \leq e^x$ for $x \in \mathbb{R}$, 
\begin{align*}
	n (\mathrm{I}) ~\lesssim~  (np)^4 (1-p)^{n-3} ~\lesssim~ \lambda^4 \exp(-\lambda/4).
\end{align*}
Moreover, from Lemma~\ref{Lemma: expectation lower bound}, we know that 
\begin{align*}
	\mE\bigl[X \sqrt{\min\{X,a\}\min\{X,b\}} \mathds{1}(X\geq 4) \bigr] \gtrsim \min\bigl\{ \lambda\sqrt{\min(\lambda,a)\min(\lambda,b)}, \lambda^4 \bigr\}.
\end{align*}
Hence, once we show 
\begin{align} \label{Eq: step 1}
	\lambda^4 \exp(-\lambda/4) \lesssim \min\bigl\{ \lambda\sqrt{\min(\lambda,a)\min(\lambda,b)}, \lambda^4 \bigr\},
\end{align}
then it implies that 
\begin{align} \label{Eq: term (I)}
	n (\mathrm{I}) \lesssim \mE\bigl[X \sqrt{\min\{X,a\}\min\{X,b\}} \mathds{1}(X\geq 4) \bigr].
\end{align}
For $\lambda < 1$, the dominating term in the upper bound on inequality~(\ref{Eq: step 1}) is $\lambda^4$ and thus inequality~(\ref{Eq: step 1}) follows easily. For $\lambda \geq 1$, by assuming $a \leq b$ without loss of generality, it is enough to show 
\begin{align*}
	& \lambda^4 \exp(-\lambda/4) \lesssim \lambda^2, \\[.5em]
	& \lambda^4 \exp(-\lambda/4) \lesssim  \lambda \sqrt{\lambda} \sqrt{a}, \quad \text{for $\lambda \leq a$}, \\[.5em]
	& \lambda^4 \exp(-\lambda/4) \lesssim \lambda \sqrt{ab}, \quad \text{for $\lambda \lesssim \min\{a,b\}$.}
\end{align*}
Starting with the first inequality, the Taylor expansion of $e^x = \sum_{n=0}^\infty \frac{x^n}{n!}$ shows that $\exp(\lambda/4) \gtrsim \lambda^2$ for all $\lambda \geq 1$. For the second inequality, again, the Taylor expansion of $e^x$ shows that $\exp(\lambda/4) \gtrsim \lambda^2 \gtrsim \lambda^{5/2}/\sqrt{a}$ for $\lambda \leq a$. The last inequality follows similarly. Hence inequality~(\ref{Eq: term (I)}) follows as desired.

\vskip 1em

\noindent \textbf{Analysis of (II).} Switching to the second term~(II), note that when $Y \geq 4$, it is always the case that $X_1 + Y \geq 4$ and $X_1' + Y \geq 4$. Therefore we can write
\begin{align*}
	& (f - f^1)^2 \mathds{1}(Y \geq 4) \\[.5em]
	=~ &  \bigl\{(X_1+ Y)\sqrt{\min\{(X_1+Y),a\}\min\{(X_1+Y),b\}} \\[.5em]
	& - (X_1'+Y)\sqrt{\min\{(X_1'+Y),a\}\min\{(X_1'+Y),b\}}\bigr\}^2 \mathds{1}(Y \geq 4). 
\end{align*}
Based on the condition that $X_1,X_1' \overset{\mathrm{i.i.d.}}{\sim} \mathrm{Binomial}(1,p)$, the expectation of $(f - f^1)^2 \mathds{1}(Y \geq 4)$ conditional on $Y$ is computed as
\begin{align*}
	\mE_{X_1,X_1'}[(f - f^1)^2 \mathds{1}(Y \geq 4)\,|\, Y]~=~ & 2p(1-p)  \bigl\{ (Y+1)\sqrt{\min\{(Y+1),a\}\min\{(Y+1),b\}} \\[.5em]
	& - Y \sqrt{\min\{Y,a\}\min\{Y,b\}} \bigr\}^2  \mathds{1}(Y \geq 4) \\[.5em]
	:= ~ & 2p(1-p) h(Y,a,b)\mathds{1}(Y \geq 4).
\end{align*}
Therefore by taking the expectation over $Y$, the second term is
\begin{align*}
	n (\mathrm{II}) ~ = ~ 2np(1-p) \mE_Y[h(Y,a,b)\mathds{1}(Y \geq 4)].
\end{align*}
For simplicity, we denote the probability mass function of $\mathrm{Binomial}(n-1,p)$ by $q_{n-1,p}(x)$ for $x = 0,1,\ldots,n-1$. Then 
\begin{align*}
	np \mE_Y[h(Y,a,b)\mathds{1}(Y \geq 4)] ~=~ np \sum_{k=4}^{n-1} h(k,a,b) q_{n-1,p}(k). 
\end{align*} 
Recall that 
\begin{align*}
	h(y,a,b) = \bigl[(y+1) \sqrt{\min\{y+1,a\}\min\{y+1,b\}} - y \sqrt{\min\{y,a\}\min\{y,b\}}\bigr]^2.
\end{align*}
and, by assuming $a \leq b$ without loss of generality and $y\geq 4$, bound this function under different cases carefully. In the following analysis, we use the fact that for $x,y \geq 4$ and some $\xi \geq 4$, 
\begin{align*}
	(y+x)^{3/2} = y^{3/2} + \frac{3\sqrt{y}}{2} x + \frac{3}{8\sqrt{\xi +y}} x^2.
\end{align*}
Therefore, it holds that $[(y+1)^{3/2} - y^{3/2}]^2 \lesssim y$ for all $y \geq 4$.	In fact, it can be numerically verified that 
\begin{align*} 
	[(y+1)^{3/2} - y^{3/2}]^2  \leq \underbrace{\frac{(5^{3/2}-4^{3/2})^2}{4}}_{:= c_0} y \quad \text{for all $y \geq 4$.}
\end{align*}
Based on this observation, we upper bound $h(y,a,b)$ as follows:
\begin{itemize}
	\item If $a \leq b \leq y$, then
	\begin{align*}
		h(y,a,b) = ab ~\leq~& \min\{y,a\} \min\{y,b\} \\[.5em]
		\leq ~ & y \sqrt{\min\{y,a\} \min\{y,b\}}.
	\end{align*}
	\item If $a \leq y < y+1 = b$, then 
	\begin{align*}
		h(y,a,b) ~=~ & [(y+1)\sqrt{a(y+1)} - y \sqrt{ay}]^2 = a [(y+1)^{3/2} -y^{3/2}]^2 \leq c_0 ay \\[.5em]
		\leq ~ & c_0 \min\{y,a\} \min\{y,b\} \\[.5em]
		\leq ~ & c_0 y \sqrt{\min\{y,a\} \min\{y,b\}}.
	\end{align*}
	\item If $a \leq y < y+1 < b$, then 
	\begin{align*}
		h(y,a,b) ~=~&  [(y+1) \sqrt{a(y+1)} - y \sqrt{ay}]^2 = a [(y+1)^{3/2} - y \sqrt{y}]^2 \leq c_0ay \\[.5em]
		\leq ~ & c_0 \min\{y,a\} \min\{y,b\} \\[.5em]
		\leq ~ & c_0 y \sqrt{\min\{y,a\} \min\{y,b\}}.
	\end{align*}
	\item If $y = a < y+1 =b$,
	\begin{align*}
		h(y,a,b) ~=~ & [(y+1) \sqrt{ab} - y \sqrt{ay}]^2 \\[.5em]
		= ~ & [(y+1) \sqrt{y(y+1)} - y^2]^2 \\[.5em]
		< ~ & [(y+1)^2 - y^2]^2  = (2y+1)^2 \leq \frac{81}{16} y^2 \quad \text{since $y \geq 4$} \\[.5em]
		= ~ & \frac{81}{16}  \min\{y,a\} \min\{y,b\} \\[.5em]
		\leq ~ &  \frac{81}{16} y \sqrt{\min\{y,a\} \min\{y,b\}}.
	\end{align*}
	\item If $y = a < y+1 < b$,
	\begin{align*}
		h(y,a,b) ~=~ & [(y+1)\sqrt{a(y+1)} - y \sqrt{ay}]^2 \\[.5em]
		= ~ & [(y+1)\sqrt{y(y+1)} - y^2 ]^2 \\[.5em]
		< ~ & [(y+1)^2 - y^2]^2  = (2y+1)^2 \leq \frac{81}{16} y^2 \quad \text{since $y \geq 4$.} \\[.5em]
		= ~ & \frac{81}{16}  \min\{y,a\} \min\{y,b\} \\[.5em]
		\leq ~ &  \frac{81}{16} y \sqrt{\min\{y,a\} \min\{y,b\}}.
	\end{align*}
	\item If $y < a = y+1 =b $,
	\begin{align*}
		h(y,a,b) ~=~& [(y+1)^2- y^2 ]^2 = (2y+1)^2 <  \frac{81}{16} y^2  \quad \text{since $y \geq 4$.} \\[.5em]
		= ~ & \frac{81}{16}  \min\{y,a\} \min\{y,b\} \\[.5em]
		\leq ~ &  \frac{81}{16} y \sqrt{\min\{y,a\} \min\{y,b\}}.
	\end{align*}
	\item If $y < y+1 < a \leq b$,
	\begin{align*}
		h(y,a,b) ~=~ & [(y+1)^2 - y^2]^2 = (2y+1)^2 < \frac{81}{16} y^2 \leq 9 \min\{b^2,y^2\}  \quad \text{since $y \geq 4$.} \\[.5em]
		= ~ & \frac{81}{16}  \min\{y,a\} \min\{y,b\} \\[.5em]
		\leq ~ &  \frac{81}{16} y \sqrt{\min\{y,a\} \min\{y,b\}}.
	\end{align*}
\end{itemize}	
In summary, for $y \geq 4$,
\begin{align*}
	h(y,a,b) \lesssim \min\{y,a\} \min\{y,b\} \lesssim y \sqrt{\min\{y,a\} \min\{y,b\}}.
\end{align*}
Using this upper bound, we have
\begin{align*}
	np \mE_Y[h(Y,a,b)\mathds{1}(Y \geq 4)] ~=~ & np \sum_{k=4}^{n-1} h(k,a,b) q_{n-1,p}(k) \\[.5em]
	\overset{\mathrm{(i)}}{\lesssim} ~ & np \sum_{k=4}^{n-1} \min\{y,a\} \min\{y,b\}  q_{n-1,p}(k) \\[.5em] 
	\overset{\mathrm{(ii)}}{\lesssim} ~ & np \sum_{k=4}^{n-1}  y \sqrt{\min\{y,a\} \min\{y,b\}} q_{n-1,p}(k). 
\end{align*}
The last inequality~(ii) yields that 
\begin{equation}
	\begin{aligned} \label{Eq: when np < 1}
		n \mathrm{(II)}~\lesssim~ & np \mE[Y\sqrt{\min\{Y,a\}\min\{Y,b\}} \mathds{1}(Y\geq 4)] \\[.5em]
		\lesssim  ~ & np \mE[X\sqrt{\min\{X,a\}\min\{X,b\}} \mathds{1}(X\geq 4)], 
	\end{aligned}
\end{equation}
where the second inequality above uses the fact that $Y\sqrt{\min\{Y,a\}\min\{Y,b\}} \mathds{1}(Y\geq 4) \leq (X_1+Y)\sqrt{\min\{(X_1+Y),a\}\min\{(X_1+Y),b\}} \mathds{1}(X_1+Y\geq 4)$. On the other hand, using inequality~(i), we can obtain
\begin{align} \label{Eq: bound on n(II)}
	n  \mathrm{(II)} ~ \lesssim \lambda \min\{ab, a \lambda, \lambda^2 + \lambda \}.  
\end{align}

\vskip 1em

\noindent \textbf{Combining pieces.} By the Efron--Stein inequality (Lemma~\ref{Lemma: Efron-Stein}), we obtained that 
\begin{align*}
	\mathrm{Var}\bigl[ X\sqrt{\min\{X,a\}\min\{X,b\}} \mathds{1}(X\geq 4) \bigr] \lesssim n \mathrm{(I)} +  n\mathrm{(II)}.
\end{align*}
By focusing on $n \mathrm{(I)}$, the bound given in (\ref{Eq: term (I)}) together with Lemma~\ref{Lemma: expectation lower bound} shows that 
\begin{align*}
	n \mathrm{(I)} \lesssim \min\{np, \sqrt{ab}\} \mE\bigl[ X\sqrt{\min\{X,a\}\min\{X,b\}} \mathds{1}(X\geq 4) \bigr].
\end{align*}
Therefore, it suffices to prove
\begin{align} \label{Eq: goal in part II}
	n\mathrm{(II)} \lesssim \min\{np, \sqrt{ab}\} \mE\bigl[ X\sqrt{\min\{X,a\}\min\{X,b\}} \mathds{1}(X\geq 4) \bigr].
\end{align}
Due to the inequality established in (\ref{Eq: when np < 1}), we know that 
\begin{align*}
	n\mathrm{(II)} \lesssim np \mE\bigl[ X\sqrt{\min\{X,a\}\min\{X,b\}} \mathds{1}(X\geq 4) \bigr] \quad \text{for all $n \geq 9$ and $p \in [0,1]$.}
\end{align*}
Hence, by assuming $\lambda = np > 1$ (when $\lambda \leq 1$, $\min\{\lambda, \sqrt{ab}\} = \lambda$),  we only need to prove that 
\begin{align*}
	n\mathrm{(II)} \lesssim  \sqrt{ab} \mE\bigl[ X\sqrt{\min\{X,a\}\min\{X,b\}} \mathds{1}(X\geq 4) \bigr],
\end{align*}
which is implied by
\begin{align*}
	\lambda \min\{ab, a \lambda, \lambda^2 + \lambda \} \lesssim \sqrt{ab} \min\bigl\{ \lambda\sqrt{\min(\lambda,a)\min(\lambda,b)}, \lambda^4 \bigr\}, 
\end{align*}
due to Lemma~\ref{Lemma: expectation lower bound} and the previous result~(\ref{Eq: bound on n(II)}). Moreover, since we assume $\lambda > 1$, the above condition is satisfied when
\begin{align} \label{Eq: sufficient condition}
	\min\{(ab)^2, a^2 \lambda^2 \} \lesssim ab \min(\lambda,a)\min(\lambda,b).
\end{align}
Now we see that 
\begin{itemize}
	\item if $\lambda \leq a \leq b$, the sufficient condition~(\ref{Eq: sufficient condition}) becomes
	\begin{align*}
		a^2\lambda^2 \lesssim ab \lambda^2  \ \Longleftrightarrow \ a \lesssim b,
	\end{align*}
	\item if $a \leq \lambda \leq b$, the sufficient condition~(\ref{Eq: sufficient condition}) becomes
	\begin{align*}
		a^2 \lambda^2 \lesssim a^2 b \lambda \ \Longleftrightarrow \ \lambda \lesssim b,
	\end{align*}
	\item if $a \leq b \leq \lambda$, the sufficient condition~(\ref{Eq: sufficient condition}) becomes
	\begin{align*}
		(ab)^2 \lesssim (ab)^2.
	\end{align*}
\end{itemize} 
Therefore we verify the sufficient condition~(\ref{Eq: sufficient condition}) and thus complete the proof of Lemma~\ref{Lemma: Variance bound Claim 22}.

\subsection{Proof of Lemma~\ref{Lemma: expectation lower bound}}
The proof of this lemma essentially follows the same lines of the proof of Claim 2.3 in \cite{canonne2018testing}. There are only two conditions that we need to check for a binomial random variable. 
\vskip 1em 

\noindent \textbf{Condition 1.} Suppose that $\lambda < 8$. For a Poisson random variable $Y \sim \mathrm{Poisson}(\lambda)$, \cite{canonne2018testing} shows that 
\begin{align} \label{Eq: Condition 1 under Pois}
	\mE\bigl[Y \sqrt{\min\{Y,a\}\min\{Y,b\}} \mathds{1}(Y\geq 4) \bigr] \geq \frac{1}{768} \min\bigl\{ \lambda\sqrt{\min(\lambda,a)\min(\lambda,b)}, \lambda^4 \bigr\}. 
\end{align}
The main step for the above result is $\mP(Y = 4) \geq \lambda^4/(4!)$. We claim that the same bound holds for $X \sim \mathrm{Binomial}(n,p)$ where $n \geq 9$ up to some constant factor. Indeed, we can write
\begin{align*}
	\mP(X=4) ~=~  & \binom{n}{4} p^4 (1-p)^{n-4} \\[.5em]
	\geq ~ & \frac{14}{729} \lambda^4 \biggl( 1 - \frac{\lambda}{n} \biggr)^{n-4} \\[.5em]
	\geq ~ & \frac{14}{729} \biggl( 1 - \frac{8}{9} \biggr)^{5} \lambda^4,
\end{align*}
where the inequalities use the fact that $n \geq 9$ and $\lambda < 8$. We then follow the same steps in \cite{canonne2018testing} and conclude that inequality~(\ref{Eq: Condition 1 under Pois}) holds when $Y$ is replaced by $X$ up to some constant factor. 

\vskip 1em

\noindent \textbf{Condition 2.} The second condition \cite{canonne2018testing} makes use of is that for $Y \sim \mathrm{Poisson}(\lambda)$ and $\lambda \geq 8$, $\mP(Y \geq \floor{\lambda/2}) \geq 1/2$. We now prove that the same inequality holds for $X \sim \mathrm{Binomial}(n,p)$. First assume that $\floor{\lambda/2} \neq 0$. Otherwise, the claim is trivial. By Chebyshev's inequality
\begin{align*}
	\mP(X < \floor{\lambda/2}) ~=~ &  \mP(-X + \lambda > \lambda - \floor{\lambda/2}) \\[.5em]
	\leq ~ & \frac{\mathrm{Var}[X]}{(\lambda - \floor{\lambda/2})^2}  \leq \frac{4\lambda}{\lambda^2} = \frac{4}{\lambda}.
\end{align*}
Now $\mP(X \geq \floor{\lambda/2}) \geq 1/2$ holds since we assume $\lambda \geq 8$. Given this probability bound, we can follow the same lines of the proof in \cite{canonne2018testing} and show that
\begin{align*}
	\mE\bigl[X \sqrt{\min\{X,a\}\min\{X,b\}} \mathds{1}(X\geq 4) \bigr] \gtrsim \min\bigl\{ \lambda\sqrt{\min(\lambda,a)\min(\lambda,b)}, \lambda^4 \bigr\}.
\end{align*}
for $\lambda \geq 8$. This completes the proof of Lemma~\ref{Lemma: expectation lower bound}.

\subsection{Proof of Lemma~\ref{Lemma: NA}} \label{Section: Proof of Lemma: NA}
Both statements follow by applying Lemma~\ref{Lemma: NA of multinomials} with specific choices of $A_1,A_2,f_1,f_2$. In particular, the first result holds by taking $A_1 = \{i\}$, $A_2 = \{j\}$, and $f_1(x) = f_2(x) = x \mathds{1}(x \geq 4)$. The second statement follows similarly by taking $A_1 = \{i\}$, $A_2 = \{j\}$, and $f_1(x) = f_2(x) = x \mathds{1}(x \geq 4) \sqrt{\min \{x, \ell_1\}\min\{x, \ell_2\}}$. 

\subsection{Proof of Lemma~\ref{Lemma: quantile}} \label{Section: Proof of Lemma: quantile}
By the definition of the quantile function and letting $V_0 := V$, $q_{1-\alpha}$ can be written as
\begin{align*}
	q_{1-\alpha} ~=~& \inf \bigg\{x \in \mathbb{R}: 1-\alpha \leq  \frac{1}{B+1} \sum_{i=0}^B \mathds{1} \big( V_i \leq x \big) \bigg\} \\[.5em]
	= ~ & V_{(k)},
\end{align*}
where $V_{(k)}$ is the $k = \ceil{(1-\alpha)(B+1)}$th order statistic of $\{V_0,V_1,\ldots,V_B\}$. Given this representation, we see that $V \leq q_{1-\alpha}$ holds if and only if
\begin{align*}
	\frac{1}{B+1} \sum_{i=0}^B \mathds{1} \bigl( V_i < V \bigr) < 1-\alpha,
\end{align*}
which is equivalent to 
\begin{align*}
	 \frac{1}{B+1} \Biggl[ \sum_{i=1}^B \mathds{1}(V_i \geq V) +1 \Biggr] > \alpha.
\end{align*}
This completes the proof of Lemma~\ref{Lemma: quantile}.

\end{document}